\documentclass{article}[a4paper]

\usepackage{amsthm} 
\usepackage{amsmath} 
\usepackage{amssymb} 
\usepackage{mathrsfs}
\usepackage{graphicx} 
\usepackage{cite} 
\usepackage{url} 
\usepackage{setspace} 
\usepackage{geometry} 
\usepackage{hyperref} 
\hypersetup{hidelinks} 
\usepackage{todonotes} 
\usepackage{color} 
\usepackage{cases}
\numberwithin{equation}{section}

\newtheorem{thm}{Theorem}[section]
\newtheorem{lem}[thm]{Lemma}
\newtheorem{defi}[thm]{Definition}

\newtheorem{pro}[thm]{Proposition}

\def\R{\mathbb{R}}
\begin{document}
	\title{Local well-posedness for   a generalized sixth-order Boussinesq equation}
	
	\author{ Long Zhong\\
		School of Mathematical Sciences,\\
		University of Electronic Science and Technology of China,\\
		Chengdu, Sichuan 610054, China \\
		\texttt{zhonglong@std.uestc.edu.cn}
		\and
		Shenghao Li \footnote{Corresponding author}\\
		School of Mathematical Sciences,\\
		University of Electronic Science and Technology of China,\\
		Chengdu, Sichuan 610054, China \\
		\texttt{lish@uestc.edu.cn}}
	\date{}
	
	\maketitle
	
\begin{abstract}
	A formally second order correct Boussinesq-type equation that describes unidirectional shallow water waves is derived,
	$$u_{tt} - u_{xx} - u_{xxxx} - u_{xxxxxx} - (u^2)_{xx} - (u^2)_{xxxx} - (uu_{xx})_{xx} - (u^3)_{xx} = 0.$$
	Such equation is analogous to original Boussinesq equation	 but with higher order approximation which may  ensure a  more accuracy description on a   long time scale. Moreover, through a rigorous derivation from Boussiensq systems, it has redeemed all the non-linear terms neglected in the sixth order Boussinesq equation (SOBE),
	$$u_{tt} - u_{xx} - u_{xxxx} - u_{xxxxxx} - (u^2)_{xx} = 0.$$ 
	 The Cauchy problem for this generalized SOBE is then considered under the Bourgain space, $X^{s,b}$, framework. The multi-linear estimates for $(u^2)_{xx}$, $(u^2)_{xxxx}$, $(uu_{xx})_{xx}$ and $(u^3)_{xx}$ are given,  the local wellposedness of the gSOBE is established for $s>\frac{1}{2}$.
\end{abstract}
	
	\section{Introduction}

In 1870s, Boussinesq introduced a uni-directional approximation model of the full Euler equation, known as the Boussinesq equation,  
\begin{equation}\label{Bou}
u_{tt} - u_{xx} - u_{xxxx}-(u^2)_{xx}=0,
\end{equation}
which helped to explain the phenomenon of solitary waves that was observed by Russell in 1834. Such equation was derived  under the assumption of shallow water waves, which requires waves to be small amplitude and long wavelength. More precisely, it needs
\[\alpha=\frac{A}{h}\ll1,\quad \beta=\frac{h^2}{l^2}\ll1,\quad S=\frac{\alpha}{\beta}\sim 1.\]
Here, $A$ denotes a typical wave amplitude, $h$ is the constant depth of water and $l$ is a typical wavelength. 

The original Boussinesq equation  has  been used in a considerable range of applications such as coast and harbor engineering, simulation of tides and tsunamis. However, it is notorious for its ill-posedness on the initial value problem (IVP). One way to correct such issue is to  establish equation from  the Euler equations with changing parts of partial derivatives for $x$ into $t$, known as the regularized Boussinesq equation,
\begin{equation}\label{rBou}
u_{tt} - u_{xx} - u_{xxtt}-(u^2)_{xx}=0.
\end{equation}
Models like these Boussiensq-type equations are known as the first order correct of the Euler equations. In fact, according to the $abcd-$system that derived from the full Euler equations by Bona, Chen and Saut \cite{bona2002boussinesq},
\begin{equation*}
\begin{cases}
\eta_t + w_x+\alpha(w\eta)_x+\beta\left(a w_{x x x}-b \eta_{x x t}\right)=0\\
w_t +\eta_x+\alpha w w_x+\beta\left(c \eta_{x x x}-d w_{x x t}\right)=0, 
\end{cases}
\end{equation*}
one can deduce both  the Boussinesq equation and the regularized Boussinesq equation as,
\[ \eta_{tt} - \eta_{xx} - \frac13\beta \eta_{xxxx}-\frac23 \alpha(\eta^2)_{xx}=0, \quad \mbox{and} \quad \eta_{tt} - \eta_{xx} - \frac13\beta \eta_{xxtt}-\frac23 \alpha(\eta^2)_{xx}=0,\]
by neglecting   terms quadratic  in $\alpha$ and $\beta$.

Another way to redeem the ill-posed issue of \eqref{Bou} is to derive a second order approximation of the Euler equations. A higher-order correct model can increase accuracy, and, just as important, are formally valid on time scales longer than ﬁrst-order 
correct models (see Section 2 of  \cite{bona2002boussinesq} for more details).
	The original  SOBE,
\begin{align} \label{1.1}
u_{tt} - u_{xx} - u_{xxxx} - u_{xxxxxx} -(u^2)_{xx} = 0,
\end{align}
was initially derived  by Christov, Maugin and Velarde \cite{christov1996well} from the full Euler equations under shallow-wave assumption with second order correct. However, part of the nonlinear terms are neglected (see \cite{christov1996well}) for simplicity of the model.  Inspired by ideas of Bona, Carvajal, Panthee, Scialom \cite{bona2018higher} and Chen \cite{chen2004boussinesq}, we derive the fully generalized SOBE\footnote[1]{
	$\Theta=\frac{1}{2}(a+d)(a+b-c-d)-(a_1-b_1+c_1-d_1)-\frac{1}{6}(a+d)-\frac{1}{3}b < 0$
}  with second order correct,
\begin{align}
\eta_{tt} - \eta_{xx} 
-\frac{1}{3} \beta\eta_{xxxx}  
-\frac{2}{3} \alpha(\eta^2)_{xx}       
+ \beta^2\Theta\eta_{xxxxxx}
- \frac{1}{2}\alpha^2(\eta^3)_{xx}               
- \alpha\beta\big[
\frac{2}{3}(\eta\eta_{xx})_{xx}  + \frac{1}{3}(\eta^2)_{xxxx}\big] = 0,
\label{1.2}
\end{align} 
from the high order  $abcd$ system introduced by Bona, Chen and Saut\cite{bona2002boussinesq},
\begin{equation*}
\left\{
\begin{aligned}
\eta_t &+w_x+\beta\left(a w_{x x x}-b \eta_{x x t}\right)+\beta^2\left(a_1 w_{x x x x x}+b_1 \eta_{x x x x t}\right) \\
&=-\alpha(\eta w)_x+\alpha\beta\big[b(\eta w)_{x x x}-(a+b-\frac{1}{3})\left(\eta w_{x x}\right)_x\big], \\
w_t &+\eta_x+\beta\left(c \eta_{x x x}-d w_{x x t}\right)+\beta^2\left(c_1 \eta_{x x x x}+d_1 w_{x x x t}\right) \\
&=-\alpha w w_x+\alpha \beta\big[(c+d) w w_{x x x}-c\left(w w_x\right)_{x x}-\left(\eta \eta_{x x}\right)_x+(c+d-1) w_x w_{x x}\big],
\end{aligned}
\right.
\end{equation*}
where $a,b,c,d,a_1,b_1,c_1,d_1$ are constants (see \cite{bona2002boussinesq,bona2018higher}). 
Rescaling (\ref{1.2}), we can obtain the gSOBE as, 
\begin{align*}
u_{tt} - u_{xx} -  u_{xxxx} - u_{xxxxxx} - (u^2)_{xx} - (u^2)_{xxxx} - (uu_{xx})_{xx} - (u^3)_{xx} = 0.
\end{align*}
	
	Thanks to the tremendous achievements on the IVPs of dispersive equations, such as the KdV equation and the Boussinesq equation in particular, the well-posedness problems of (\ref{1.1}) has been well-studied for both IVP and initial boundary value problem (IBVP) during the past decades. Esfahani, Farah, and Wang \cite{esfahani2012global} first establish the well-posedness of
	\begin{equation*}
		\begin{cases}
			u_{tt} - u_{xx} + k u_{xxxx} - u_{xxxxxx} - (u^2)_{xx} = 0, \quad x \in \mathbb{R},t \geq 0, \\
			u(x,0) = \varphi(x), \quad u_t(x,0) = \psi''(x),
		\end{cases}
	\end{equation*}
	in $H^{s}(\mathbb{R})$ for initial data $(\varphi, \psi) \in H^s(\mathbb{R}) \times H^{s-1}(\mathbb{R})$ with $s\geq0$ and $k = \pm 1$ by employing the Strichartz estimate. Next, Esfahani and Farah \cite{esfahani2012local} improve their result to $s > -\frac{1}{2}$ using a related Bourgain space framework. Finally, Esfahani and Wang \cite{esfahani2014bilinear} extend the conclusion to $s > -\frac{3}{4}$ based on the $[k; Z]$-multiplier norm method introduced by Tao \cite{tao2001multilinear}. The related IBVPs are also studied in \cite{li2018nonhomogeneous,li2019wellposedness,li2020nonhomogeneous,li2022lower} with related space accordingly. 
	
	
	In this article, we study the Cauchy problem of the gSOBE
	\begin{equation}
		\begin{cases}
			u_{tt} - u_{xx} +k  u_{xxxx} - u_{xxxxxx} - (u^2)_{xx} - (u^2)_{xxxx} - (uu_{xx})_{xx} - (u^3)_{xx} = 0, \quad x \in \mathbb{R},t \geq 0, \\
			u(x,0) = \varphi(x), \quad u_t(x,0) = \psi''(x),
		\end{cases}
		\label{1.6}
	\end{equation}
	with $k = \pm 1$. Before presenting the main theorem, we denote $\langle x \rangle = \sqrt{1+x^2}$ and introduce the corresponding Bougain-type space, $X^{s,b}$, for the sixth-order Boussinesq equation.
\begin{defi}\label{defn}
For $s, b \in \mathbb{R}$, $X^{s,b}$ denotes the completion of the Schwartz class $S(\mathbb{R}^2)$ with
\[
\|u\|_{X^{s,b}} = \|\langle \xi \rangle ^s \langle |\tau| - \phi(\xi) \rangle^b \hat u(\xi,\tau)\|_{L^2_{\xi,\tau}(\mathbb{R}^2)}
\]
where $\phi(\xi) = \sqrt{\xi^6 + k\xi^4 +\xi^2}$ and $\hat u$ denotes the Fourier transform on both time and space of $u$. 
\end{defi} 
\noindent	We then set $X^{s,b}_T := X^{s,b}|_{\mathbb{R}\times(0,T)}$ with the quotient norm,
	\[
	\|u\|_{X^{s,b}(\mathbb{R}\times(0,T))} := \inf_{w \in X^{s,b}} \{\|w\|_{X^{s,b}} : w(x,t) = u(x,t) \quad \text{on}\quad \mathbb{R}\times(0,T)\}.
	\]
	 The main result of this paper is as follows.
\begin{thm}\label{main}
 Let $s>\frac{1}{2}$ be given, there exists $b = b(s) \in (\frac{1}{2}, 1)$, we can find $T = T(s, b)>0$ such that if
$$(\varphi, \psi) \in H^s(\mathbb{R}) \times H^{s-1}(\mathbb{R}),$$
then the IVP(\ref{1.6}) admits a unique solution
$u \in X^{s,b}_T,$
and the corresponding solution map is real analytic.
\end{thm}

	Recalling from previous results on the linear SOBE, 
	\begin{equation}
		\begin{cases}
			u_{t t}- u_{x x}+k u_{x x x x}-u_{x x x x x x}=0, \quad x \in \mathbb{R}, t>0, \\
			u(x, 0)=\varphi(x), u_t(x, 0)=\psi^{\prime \prime}(x) ,
		\end{cases}
		\label{ivp}
	\end{equation}
	we write the solution $u =\left[V_1\left(\varphi\right)\right](x, t)+\left[V_2\left(\psi\right)\right](x, t)$ with
	$$
	\left[V_1\left(\varphi\right)\right](x, t):=\frac{1}{2} \int_{\mathbb{R}}\left(e^{i(t \phi(\xi)+x \xi)}+e^{i(-t \phi(\xi)+x \xi)}\right) \widehat{\varphi}(\xi) d \xi,
	$$
	and
	$$
	\left[V_2\left(\psi\right)\right](x, t):=\frac{1}{2 i} \int_{\mathbb{R}}\left(e^{i(t \phi(\xi)+x \xi)}-e^{i(-t \phi(\xi)+x \xi)}\right) \frac{\xi^2 \widehat{\psi}(\xi)}{\phi(\xi)} d \xi .
	$$
	According to the Duhamel's principal, the IVP for the forced linear equation,
	$$
	\left\{\begin{array}{l}
		u_{t t}- u_{x x}+ku_{x x x x}-u_{x x x x x x}=f(x, t), \quad x \in \mathbb{R}, t>0, \\
		u(x, 0)=0, u_t(x, 0)=0
	\end{array}\right.
	$$
	has its solution $u$ in the form
	$$
	u(x, t)=\int_0^t\left[V_2(f)\right]\left(x, t-t^{\prime}\right) d t^{\prime} .
	$$
	Let us denote $\eta(t)$ to be a cut-off function such that $\eta \in C_0^{\infty}(\mathbb{R})$ with $\eta(t)=1$ on $(-1,1)$ and supp $\eta \in(-2,2)$. The estimates below come from \cite{farah2009local,esfahani2012local}.
\begin{lem}\label{l1}
	For any $s, b \in \mathbb{R}$,  the solution $u$ of the IVP (\ref{ivp}) satisfies
	$$
	\begin{aligned}
	\|\eta(t) u(x, t)\|_{X^{s, b}}  \lesssim \left\|\varphi\right\|_{H^s(\mathbb{R})}+\left\|\psi\right\|_{H^{s-1}(\mathbb{R})}.
	\end{aligned}
	$$
\end{lem}
 \begin{lem}\label{l2}
 Let $-\frac{1}{2}<b^{\prime} \leq 0 \leq b \leq b^{\prime}+1$ and $0<T \leq 1$ then
 $$
 \left\|\eta\left(\frac{t}{T}\right) \int_0^t\left[V_2(f)\right]\left(x, t-t^{\prime}\right) d t^{\prime}\right\|_{X^{s, b}} \leq  T^{1+b^{\prime}-b}\left\|\left(\frac{\widehat{f}(\xi, \tau)}{\phi(\xi)}\right)^{\vee}\right\|_{X^{s, b^{\prime}}},
 $$
 where "$\vee$" denote the inverse Fourier transform in both time and space.
 \end{lem}

	Through a standard contraction mapping principal technique, Theorem \ref{main} can be obtained from Lemmas \ref{l1}-\ref{l2}(see \cite{farah2009local, esfahani2012local} for details) once related   multi-linear estimates are established. The method of multi-linear estimates on $X^{s,b}$ space is originally introduced by Kenig-Ponce-Vega \cite{kenig1996bilinear} for the KdV equation. Recently, based on Tao's $[k; Z]$-multiplier method \cite{tao2001multilinear}, Li, Yang and Zhang \cite{li2023lower,yang2022local, li2023nonhomogeneous, li2018nonhomogeneous,li2019wellposedness,li2020nonhomogeneous,li2022lower} have refined a method for bi-linear estimates that can fulfill for different values of $b$ in $X^{s,b}$ space (including $b < \frac{1}{2}$, $b = \frac{1}{2}$, $b > \frac{1}{2}$). Such technique on the $X^{s,b}-$ space for variety choices of $b$ is essential for both IVP and IBVP(see \cite{li2023lower,yang2022local, li2023nonhomogeneous, li2018nonhomogeneous,li2019wellposedness,li2020nonhomogeneous,li2022lower}). We will adapt this approach to handle bilinear estimates on terms, $(u^2)_{xx}$, $(u^2)_{xxxx}$, and $(uu_{xx})_{xx}$. In addition, inspired by the multi-linear estimate of the coupled KdV system \cite{carvajal2019sharp}, we will generalize such refined method and help to establish the trilinear estimate for term $(u^3)_{xx}$. 

\begin{pro}\label{bi01}
 For any $s > -\frac{3}{4}$ , there exists a $b_0 > \frac{1}{2}$ such that for any $b \in (\frac{1}{2},b_0],$

\begin{align}
\left\|\left(\frac{\widehat{(u^2)_{xx}}(\tau,\xi)}{\phi(\xi)}\right)^{\vee}\right\|_{X^{s, b-1}} & \lesssim\|u\|^2_{X^{s, b}} . \label{bi-1}\\
	\left\|\left(\frac{\widehat{(u^2)_{xxxx}}(\tau,\xi)}{\phi(\xi)}\right)^{\vee}\right\|_{X^{s, b-1}} & \lesssim\|u\|^2_{X^{s, b}} . \label{bi-2}
	\end{align}
\end{pro}

\begin{pro}\label{bi03}
For any $s > \frac{1}{2}$ , there exists a $b_0 > \frac{1}{2}$ such that for any $b \in (\frac{1}{2},b_0],$

\begin{align}
\left\|\left(\frac{\widehat{(uu_{xx})_{xx}}(\tau,\xi)}{\phi(\xi)}\right)^{\vee}\right\|_{X^{s, b-1}} & \lesssim\|u\|^2_{X^{s, b}}. \label{bi-3}\\
\left\|\left(\frac{\widehat{(u^3)_{xx}}(\tau,\xi)}{\phi(\xi)}\right)^{\vee}\right\|_{X^{s, b-1}} & \lesssim \|u\|^3_{X^{s, b}}.
\label{tri-1} 
	\end{align}
\end{pro}

	The paper is organized as follows: We first derive the gSOBE in Section 2. Next, desired multi-linear estimates will be established in Section 3.
	\section{Derivation of the gSOBE}
	Starting from bi-directional model of the Boussinesq system\footnote[2]{$\alpha,\beta, a,b,c,d,a_1,b_1,c_1,d_1$ are constants.}
	\begin{align}{}
		\eta_t& + w_x+\beta\left(a w_{x x x}-b \eta_{x x t}\right)+\beta^2\left(a_1 w_{x x x x x}+b_1 \eta_{x x x x t}\right)+\alpha(\eta w)_x \label{2.1a} \\
		&-\alpha\beta\big[b(\eta w)_{x x x}-(a+b-\frac{1}{3})\left(\eta w_{x x}\right)_x\big] = O(\alpha^3, \beta^3,\alpha^2\beta, \alpha\beta^2), \nonumber\\
		w_t& +\eta_x+\beta\left(c \eta_{x x x}-d w_{x x t}\right)+\beta^2\left(c_1 \eta_{x x x x}+d_1 w_{x x x t}\right)+\alpha w w_x \label{2.1b}\\
		&-\alpha \beta\big[(c+d) w w_{x x x}-c\left(w w_x\right)_{x x}-\left(\eta \eta_{x x}\right)_x+(c+d-1) w_x w_{x x}\big] = O(\alpha^3, \beta^3,\alpha^2\beta, \alpha\beta^2), \nonumber
	\end{align}
	which was introduced in \cite{bona2002boussinesq}, we derive the gSOBE. Note that $a+b+c+d=\frac{1}{3}$ and $\alpha, \beta$ are small parameters with $\alpha, \beta \ll 1, \alpha/\beta\sim 1.$
	Differentiating (\ref{2.1a}) and (\ref{2.1b}) with respect to $t$ and $x$ accordingly, and taking their difference, one has
	\begin{align}
		\eta_{tt} - \eta_{xx} + \beta \big[(a+d) w_{xxxt}-b\eta_{xxtt}-c\eta_{xxxx}\big]+\beta^2\big[(a_1-d_1) w_{xxxxxt}+b_1\eta_{xxxxtt}-c_1\eta_{xxxxxx}\big] \nonumber\\
		+\alpha [(\eta w)_{xt}-( w w_x)_x]+\alpha\beta\big[b(\eta w)_{xxxt}-(a+b-\frac{1}{3})(\eta w_{xx})_{xt}-(c+d)( w w_{xxx})_x \nonumber\\
		+c( w w_x)_{xxx}+(\eta\eta_{xx})_{xx}-(c+d-1)( w_x w_{xx})_x\big] = O(\alpha^3, \beta^3,\alpha^2\beta, \alpha\beta^2).
		\label{2.2}
	\end{align}
	Next, we intend to derive a single equation with respect to $\eta$, which is similar to the SOBE (\ref{1.1}). This requires using linear and nonlinear combination of $\partial_x^n\eta$, $n = 0,1,2,3$, to approximate $ w$ and $\eta_t$.
	Since $\alpha,\beta$ are small and $\alpha \sim \beta$, we first consider the system (\ref{2.1a})-(\ref{2.1b}) in $O(1)$ with respect to $\alpha$ and $\beta$, which reads like the one-dimensional wave equation as,
	\begin{align*}
		\eta_t+w_x=0, \quad w_t+\eta_x=0.
	\end{align*}
	 Given initial condition $\eta(x, 0)=f(x), w(x, 0)=g(x)$\footnote[3]{$f(x)$ and $g(x)$ are the initial disturbances of the surface and the horizontal velocity, respectively.},
	the solution can be obtained as,
	$$
	\left\{\begin{aligned}
		\eta(x, t) & =\frac{1}{2}\big[f(x+t)-g(x+t)\big]+\frac{1}{2}\big[f(x-t)+g(x-t)\big], \\
		w(x, t) & =\frac{1}{2}\big[g(x+t)-f(x+t)\big]+\frac{1}{2}\big[g(x-t)+f(x-t)\big] .
	\end{aligned}\right.
	$$
	Considering the solution is unidirectional, without loss of generality, we assume the wave is only \textbf{moving to the right}, which leads to $f=g$ and $\eta =  w = f(x-t)$. Thus, we have the desired first order approximations  in the Boussinesq system,
	\begin{align}
		\eta_t = -\eta_x + O(\alpha, \beta), \quad  w = \eta + O(\alpha, \beta), \quad  \text{as} \quad \alpha, \beta \to 0.		\label{2.3}
	\end{align}
	Furthermore, we consider a second order approximation based on (\ref{2.1a})-(\ref{2.1b}) writing as, 
	\begin{align}
		 w = \eta + \alpha A +\beta B + O(\alpha^2, \beta^2,\alpha\beta), \quad \text{as} \quad \alpha, \beta \to 0 \label{2.4}
	\end{align}
	where $A$ and $B$ are to be determined. Substituting (\ref{2.4}) into (\ref{2.1a})-(\ref{2.1b}), one can obtain
	\begin{align}
		 w = \eta - \frac{1}{4} \alpha \eta^2 + \frac{1}{2}\beta (c+d-a-b)\eta_{xx}+ O(\alpha^2, \beta^2,\alpha\beta), \quad \text{as} \quad \alpha, \beta \to 0. \label{2.5}
	\end{align}
	Then, plugging (\ref{2.5}) into (\ref{2.1a}) will leads to, 
	\begin{align}
		\eta_t= - \eta_x-\frac{3}{2}\alpha\eta\eta_x -\frac{1}{6}\beta\eta_{xxx}+O(\alpha^2, \beta^2,\alpha\beta), \quad \text{as} \quad \alpha, \beta \to 0. \label{2.6} 
	\end{align}
	For more detail, interested reader can refer to   Section 2 in \cite{bona2002boussinesq}.
For simplicity in notations, we will neglect terms with $O(\alpha^3, \beta^3,\alpha^2\beta, \alpha\beta^2)$. Substituting (\ref{2.3}) and (\ref{2.5}) into  (\ref{2.2}) leads to
	\begin{align}
		\eta_{tt} - \eta_{xx} 
		+ \beta\big[(a+d)\eta_{xxxt}-b\eta_{xxtt}-c\eta_{xxxx}\big]  
		+ \alpha\big[(\eta^2)_{xt} - \frac{1}{2}(\eta^2)_{xx}\big]       
		+ \beta^2\Theta_1\eta_{xxxxxx}
		+ \alpha^2(\frac{1}{2}\eta^3)_{xx}   \nonumber  \\          
		+ \alpha\beta\big[(c+d)(\eta\eta_{xxx})_x + (c+d-1)(\eta_x\eta_{xx})_{x} 
		- (c+d+\frac{2}{3})(\eta\eta_{xx})_{xx}  + \frac{1}{4}(a+4b-2c+d)(\eta^2)_{xxxx}\big] = 0
		\label{2.7}
	\end{align}
	where $\Theta_1=\frac{1}{2}(a+d)(a+b-c-d)-(a_1-b_1+c_1-d_1)$.
	In addition, according to (\ref{2.6}), one has
	\begin{equation}
		\left\{
		\begin{aligned}
			(\eta^2)_{xt} &= -(\eta^2)_{xx}-3\alpha(\eta^2\eta_x)_x-\frac{1}{3}\beta(\eta\eta_{x x x})_x ,\\
			\eta_{x x x t} &= -(\eta_x+\frac{3}{2}\alpha\eta\eta_x+\frac{1}{6}\beta\eta_{xxx})_{xxx},\\
			\eta_{xxtt} &= \eta_{x x x x} +3\alpha(\eta\eta_x)_{xxx}+\frac{1}{3}\beta\eta_{xxxxx}.
		\end{aligned}
		\right.
		\label{2.8}
	\end{equation}
	Then, substituting (\ref{2.8}) into (\ref{2.7})  follows that
	\begin{align}
		\eta_{tt} - \eta_{xx} 
		-\frac{1}{3} \beta\eta_{xxxx}  
		-\frac{2}{3} \alpha(\eta^2)_{xx}       
		+ \beta^2\Theta\eta_{xxxxxx}
		- \frac{1}{2}\alpha^2(\eta^3)_{xx}               
		+ \alpha\beta[(c+d-\frac{1}{3})(\eta\eta_{xxx})_x + \nonumber\\ (c+d-1)(\eta_x\eta_{xx})_{x} 
		- (c+d+\frac{2}{3})(\eta\eta_{xx})_{xx}  + \frac{1}{4}(-2a-2b-2c-2d)(\eta^2)_{xxxx}] = 0,
		\label{2.9}
	\end{align}
	with $\Theta=\frac{1}{2}(a+d)(a+b-c-d)-(a_1-b_1+c_1-d_1)-\frac{1}{6}(a+d)-\frac{1}{3}b < 0$. Moreover, based on the chain rule for differentiation, we have
	\begin{align*}
		\eta_x\eta_{xx}+\eta\eta_{xxx}=(\eta\eta_{xx})_{x},\quad \eta\eta_{xxx} = \frac{1}{2}(\eta^2)_{xxx}-\frac{3}{2}(\eta^2_x)_x. 
	\end{align*} 
	Thus, (\ref{2.9}) can be rewritten as
	\begin{align}
		\eta_{tt} - \eta_{xx} 
		-\frac{1}{3} \beta\eta_{xxxx}  
		-\frac{2}{3} \alpha(\eta^2)_{xx}       
		+ \beta^2\Theta\eta_{xxxxxx}
		- \frac{1}{2}\alpha^2(\eta^3)_{xx}               
		- \alpha\beta\big[
		 \frac{2}{3}(\eta\eta_{xx})_{xx}  + \frac{1}{3}(\eta^2)_{xxxx}\big] = 0. \label{2.10}
	\end{align}
	Therefore, we obtain the desired fully gSOBE by rescaling (\ref{2.10}),
	\begin{align*}
		u_{tt} - u_{xx} -  u_{xxxx} - u_{xxxxxx} - (u^2)_{xx} - (u^2)_{xxxx} - (uu_{xx})_{xx} - (u^3)_{xx} = 0.
	\end{align*}
	
	\section{Multi-linear estimate}
	In this section, we establish the desired multi-linear estimates, which   equivalent to the following propositions. 
	\begin{pro}\label{bi0}
		For any $s > -\frac{3}{4}$, there exists a $b_0 > \frac{1}{2}$ such that for any $b \in (\frac{1}{2},b_0],$
		\begin{align}
		\left\|\left(\frac{\xi^2\widehat{u v}}{\phi(\xi)}\right)^{\vee}\right\|_{X^{s, b-1}} & \lesssim\|u\|_{X^{s, b}}\|v\|_{X^{s, b}},\label{bi11}\\ 
	\left\|\left(\frac{\xi^4\widehat{u v}}{\phi(\xi)}\right)^{\vee}\right\|_{X^{s, b-1}} & \lesssim\|u\|_{X^{s, b}}\|v\|_{X^{s, b}}. \label{bi12}
		\end{align}
\end{pro}

	\begin{pro}
For any $s > \frac{1}{2}$, there exists a $b_0 > \frac{1}{2}$ such that for any $b \in (\frac{1}{2},b_0],$

\begin{align}
\left\|\left(\frac{\xi^2\widehat{u v_{xx}}}{\phi(\xi)}\right)^{\vee}\right\|_{X^{s, b-1}} & \lesssim\|u\|_{X^{s, b}}\|v\|_{X^{s, b}}. \label{bi13}\\
	\left\|\left(\frac{\xi^2\widehat{u v w}}{\phi(\xi)}\right)^{\vee}\right\|_{X^{s, b-1}} & \lesssim\|u\|_{X^{s, b}}\|v\|_{X^{s, b}}\|w\|_{X^{s, b}}. \label{tri1}
	\end{align}
\end{pro}

	 The idea of proofs is inherited from studies for well-posedness problem of the  Schrödinger equation, KdV-KdV, and SOBE \cite{carvajal2019sharp,li2023lower,yang2022local, li2023nonhomogeneous, li2018nonhomogeneous,li2019wellposedness,li2020nonhomogeneous,li2022lower}.
	\subsection{Preliminary}
	In order to achieve a better cancellation in   multi-linear estimates, we introduce an equivalence relation in the related Bourgain type space,
	\[
	\|u\|_{X^{s,b}(\mathbb{R}^2)} \sim \|\langle \xi \rangle ^s \langle |\tau| - |\xi|^3 - \frac{k}{2}|\xi| \rangle^b \hat u(\xi,\tau)\|_{L^2_{\xi,\tau}(\mathbb{R}^2)} ,
	\]
	which is inherited from Farah \cite{esfahani2012local,farah2009local}. For simplicity, we only present the proof for $k=-1$ since the case when $k = 1$ can be addressed similarly. In addition,  some elementary auxiliary lemmas will be adapted, which can be found in \cite{erdougan2016regularity,li2022lower,yang2022local}.
\begin{lem}\label{lem1}
	Let $\rho \geq \gamma \geq 0$ and $\rho + \gamma > 1$, then
	\[
	\int_{\mathbb{R}} \frac{1}{\langle x-c_1 \rangle^{\rho}\langle x-c_2\rangle^{\gamma}} dx \lesssim \langle c_1 - c_2 \rangle ^{-\gamma}\phi_{\rho}(c_1-c_2),\quad 
	\mbox{where} \quad 
	\phi_{\rho}(a) = 
	\begin{cases}
	1, & \text{$\rho > 1$}, \\
	log(1+\langle a\rangle), & \text{$\rho  = 1$}, \\
	\langle a \rangle ^{1-\rho}, & \text{$\rho < 1 $.}
	\end{cases}
	\]
\end{lem}
\begin{lem}\label{lem2}
	If $p > 1$, then for any $c_i \in \mathbb{R}$, $0 \leq i \leq 2$, with $c_2 \neq 0$, we have
	\[
	\int_{\mathbb{R}} \frac{1}{\langle c_2x^2 + c_1x +c_0 \rangle^{p}} dx \lesssim \frac{1}{|c_2|^{1/2}\langle c_0 -\frac{c_1^2} {4c_2} \rangle ^{1/2}}.
	\]
\end{lem}
\begin{lem}\label{lem3}
	If $p > 1/2$, then for any $c_i \in \mathbb{R}$, $0 \leq i \leq 2$, with $c_2 \neq 0$, we have
	\[
	\int_{\mathbb{R}} \frac{1}{\langle c_2x^2 + c_1x +c_0 \rangle^{p}} dx \lesssim \frac{1}{|c_2|^{1/2}}.
	\]
	Similarly, if $p > \frac{1}{3}$, then for any $c_i \in \mathbb{R}$, $0 \leq i \leq 3$, with $c_3 \neq 0$, we have
	\[
	\int_{\mathbb{R}} \frac{1}{\langle c_3x^3+c_2x^2 + c_1x +c_0 \rangle^{p}} dx \lesssim \frac{1}{|c_3|^{1/3}}.
	\]
\end{lem}
To establish multi-linear estimates, we will need to transfer them into some weighted convolution of $L^2$ functions.   For convenience of the notation, for $\vec{\xi} =\left(\xi_1,..., \xi_k\right)$ and $ \vec{\tau} =\left(\tau_1,..., \tau_k\right)$, we introduce the set $A_k$ and denote
\begin{align*}
A_k(\R):=\left\{(\vec{\xi}, \vec{\tau}) \in \R^{2k}: \sum_{i=1}^k \xi_i=\sum_{i=1}^k \tau_i=0\right\} , \quad \mbox{with } k =   3, 4.
\end{align*}

	\subsection{Proof of Proposition \ref{bi0}}

	The proof of estimate \eqref{bi11} is similar to the one for  \eqref{bi12}, therefore we only present the proof for  \eqref{bi12}. Set $p=-s$, then according to $\phi(\xi)= |\xi|\sqrt{\xi^4+\xi^2+1} \sim|\xi|\langle\xi\rangle^2$ and Definition \ref{defn}, one has
	\begin{align*}
		\left\|\left(\frac{\xi^4\widehat{u v}}{\phi(\xi)}\right)^{\vee}\right\|_{X^{s, b-1}}
		\sim \left\|\langle \xi \rangle ^s \langle|\tau| - |\xi|^3 +\frac{1}{2}|\xi| \rangle^{b - 1} \frac{\xi^4\widehat{u v}}{|\xi|\langle\xi\rangle^2}\right\|_{L^2_{\xi,\tau}(\mathbb{R}^2)} .
	\end{align*}
	Thus, based on duality and Plancherel theorem, \eqref{bi12} is equivalent to
	\begin{equation} 
		\int_{A_3(\mathbb{R}) } \left| \frac{|\xi_3|^3\left\langle\xi_1\right\rangle^p\left\langle\xi_2\right\rangle^p h_1(\xi_1,\tau_1)h_2(\xi_2,\tau_2)h_3\left(\xi_3, \tau_3\right)}{\left\langle\xi_3\right\rangle^{2+p}\left\langle L_1\right\rangle^b\left\langle L_2\right\rangle^b\left\langle L_3\right\rangle^{1-b}}  \right|\lesssim \prod_{j=1}^3\left\|h_j\right\|_{L_{\xi, \tau}^2},  \quad  h_1,h_2,h_3 \in L^2 \label{3.4}
	\end{equation}
	where $h_1 =\langle\xi\rangle^s\left\langle |\tau| - |\xi|^3 +\frac{1}{2}|\xi|\right\rangle^b \hat{u}(\xi, \tau)$, $h_2 =\langle\xi\rangle^s\left\langle |\tau| - |\xi|^3 +\frac{1}{2}|\xi|\right\rangle^b \hat{v}(\xi, \tau)$ and
	$$
	L_i:=|\tau_i| - |\xi_i|^3 +\frac{1}{2}|\xi_i|, \quad i=1,2,3.
	$$
	 In addition, we also introduce the resonance function $$H(\xi_1,\xi_2,\xi_3)=L_1+L_2+L_3,$$ 
	 which is the key to establish multi-linear estimates. Since $\left\langle\xi_1\right\rangle\left\langle\xi_2\right\rangle /\left\langle\xi_3\right\rangle \geq 1$, we only need to consider the case when $p$ is close to $\frac{3}{4}$. Without loss of generality, we assume $\frac{9}{16}<p<\frac{3}{4}$ and define $b_0 = \frac{3}{4} - \frac{1}{3}p < \frac{9}{16}$ to help to justify the proof for $\frac{1}{2}<b \leq b_0$.
	In order to establish the estimate (\ref{bi12}), we will examine the six potential cases regarding the signs of $\tau_1, \tau_2$, and $\tau_3$.
	\begin{align*}
		&\textbf{A} := \{\tau_1 \geq 0, \quad \tau_2 \geq 0, \quad \tau_3 \leq 0\}, 
		&\textbf{B} := \{\tau_1 \geq 0, \quad \tau_2 \leq 0, \quad \tau_3 \geq 0\}, \\
		&\textbf{C} := \{\tau_1 \geq 0, \quad \tau_2 \leq 0, \quad \tau_3 \leq 0\}, 
		&\textbf{D} := \{\tau_1 \leq 0, \quad \tau_2 \geq 0, \quad \tau_3 \geq 0\}, \\
		&\textbf{E} := \{\tau_1 \leq 0, \quad \tau_2 \geq 0, \quad \tau_3 \leq 0\}, 
		&\textbf{F} := \{\tau_1 \leq 0, \quad \tau_2 \leq 0, \quad \tau_3 \geq 0\}. 
	\end{align*}
	Note that cases where $\tau_1, \tau_2, \tau_3 \geq 0$ and $\tau_1, \tau_2, \tau_3 \leq 0$ are impossible due to the constraint $\tau_1 + \tau_2 + \tau_3 = 0$. To simplify the analysis, we observe that conditions $\mathbf{A}, \mathbf{B}, \mathbf{C}$ are equivalent to $\mathbf{F}, \mathbf{E}, \mathbf{D}$ respectively by considering the transformation $\left(\tau_1, \tau_2, \tau_3\right) \rightarrow -\left(\tau_1, \tau_2, \tau_3\right)$ and the fact that $L^2$ norms are invariant under reflections. Moreover, due to the symmetric structure between $\tau_1$ and $\tau_2$, cases $\mathbf{B}$ and $\mathbf{D}$ are also equivalent.
	Therefore, it suffices to prove the estimate in cases $\mathbf{A}$ and $\mathbf{B}$. For different conditions of $\left(\tau_1, \tau_2, \tau_3\right)$, we split the space of $\left(\xi_1, \xi_2, \xi_3\right) \in \mathbb{R}^3$ as follows:
	\begin{align*}
		&\textbf{(a)} := \{\xi_1 \geq 0, \quad \xi_2 \geq 0, \quad \xi_3 \leq 0\}, 
		&\textbf{(b)} := \{\xi_1 \geq 0, \quad \xi_2 \leq 0, \quad \xi_3 \leq 0\}, \\
		&\textbf{(c)} := \{\xi_1 \leq 0, \quad \xi_2 \geq 0, \quad \xi_3 \leq 0\}, 
		&\textbf{(d)} := \{\xi_1 \leq 0, \quad \xi_2 \leq 0, \quad \xi_3 \geq 0\}, \\
		&\textbf{(e)} := \{\xi_1 \leq 0, \quad \xi_2 \geq 0, \quad \xi_3 \geq 0\}, 
		&\textbf{(f)} := \{\xi_1 \geq 0, \quad \xi_2 \leq 0, \quad \xi_3 \geq 0\}. 
	\end{align*}
	It then suffices to consider cases \textbf{(a)}, \textbf{(b)} and \textbf{(c)}. Firstly, for case \textbf{A}, we can write
	$$
	L_1=\tau_1-\left|\xi_1\right|^3+\frac{1}{2}\left|\xi_1\right| ; \quad L_2=\tau_2-\left|\xi_2\right|^3+\frac{1}{2}\left|\xi_2\right| ; \quad L_3:=\tau_3+\left|\xi_3\right|^3-\frac{1}{2}\left|\xi_3\right| .
	$$
	\textbf{Case A1.} If $\left|\xi_1\right| \leq 1$, it indicates  that$\langle\xi_1\rangle\langle\xi_2\rangle \langle\xi_3\rangle \lesssim 1$. Thus, (\ref{3.4}) can be simplified to
		\begin{align}
			\int_{A_3} \frac{\left|\xi_3\right| \prod_{j=1}^3\left|h_j\left(\xi_j, \tau_j\right)\right|}{\left\langle L_1\right\rangle^{b}\left\langle L_2\right\rangle^{b}\left\langle L_3\right\rangle^{1-b}} \lesssim \prod_{j=1}^3\left\|h_j\right\|_{L_{\xi, \tau}^2}.
			\label{3.5}
		\end{align}
	\begin{itemize}
		\item \textbf{Case A1.1.}  If $\left|\xi_2\right| \leq 2$, one has $\left|\xi_3\right| \leq 3$.	Hence,
		\begin{align*}
			\text { LHS of }(\ref{3.5})  \lesssim \iint \frac{\left|h_3\right|}{\left\langle L_3\right\rangle^{1-b}}\left(\iint \frac{d \xi_1 d \tau_1}{\left\langle L_1\right\rangle^{2 b}\left\langle L_2\right\rangle^{2b}}\right)^{\frac{1}{2}}\left(\iint h_1^2 h_2^2 d \xi_1 d \tau_1\right)^{\frac{1}{2}} d \xi_3 d \tau_3 .
		\end{align*}

		According to Lemma \ref{lem1}, for $\frac{1}{2}<b<1$ and $|\xi_1|\leq 1$, one has
		$$
		\sup _{\xi_3, \tau_3} \frac{1}{\left\langle L_3\right\rangle^{2-2 b}} \int \frac{d \xi_1}{\left\langle L_1+L_2\right\rangle^{2b}} \lesssim 1 .
		$$
		
		Hence, (\ref{3.5}) is achieved through Cauchy-Schwarz inequality,
		\begin{align*}
			\text { LHS of }(\ref{3.5}) & \lesssim \iint\left|h_3\right|\left(\iint h_1^2 h_2^2 d \xi_1 d \tau_1\right)^{\frac{1}{2}} d \xi_3 d \tau_3 
			\lesssim \prod_{j=1}^3\left\|h_j\right\|_{L_{\xi,\tau}^2}.
		\end{align*}
		
		\item \textbf{Case A1.2.} If $\left|\xi_2\right|>2$, similar to the argument in \textbf{Case A1.1}, it suffices to bound  
		\begin{align}
			\sup _{\xi_3, \tau_3} \frac{\left|\xi_3\right|^2}{\left\langle L_3\right\rangle^{2-2b}} \int \frac{d \xi_1}{\left\langle L_1+L_2\right\rangle^{2 b}}.
			\label{3.6}
		\end{align}
		In addition, we notice that $\left|\xi_1\right| \leq 1$ and $\left|\xi_2\right|>2$ will lead to $\left|\xi_2\right| \sim\left|\xi_3\right|$.
		\begin{itemize}
			\item \textbf{Case A1.2.1} If $\left(\xi_1, \xi_2, \xi_3\right)$ lies in \textbf{(a)}, that is, $\xi_1 \geq 0$, $\xi_2 \geq 0$, $\xi_3 \leq 0.$ Then, one has
			$$
			L_1+L_2=3 \xi_3 \xi_1^2+3 \xi_3^2 \xi_1-L_3, \quad \text { with } \quad L_3=\tau_3-\xi_3^3+\frac{1}{2} \xi_3.
			$$
			
			For $b \leq b_0<\frac{3}{4}$, it yields $2-2b >\frac{1}{2}$. Then, according to Lemma \ref{lem2}, one has
			$$
			(\ref{3.6}) \lesssim \frac{\left|\xi_3\right|^2}{\left\langle L_3\right\rangle^{2-2b}} \frac{1}{\left|\xi_3\right|^{\frac{1}{2}}}\left\langle L_3+\frac{3}{4} \xi_3^3\right\rangle^{-\frac{1}{2}} \lesssim \frac{\left|\xi_3\right|^{\frac{3}{2}}}{\left\langle L_3\right\rangle^{2-2b}\left\langle L_3+\frac{3}{4} \xi_3^3\right\rangle^{\frac{1}{2}}} \lesssim 1,
			$$
			since $\left\langle L_3\right\rangle^{2-2b}\left\langle L_3+\xi_3^3\right\rangle^{\frac{1}{2}} \gtrsim\left\langle L_3\right\rangle^{\frac{1}{2}}\left\langle L_3+\xi_3^3\right\rangle^{\frac{1}{2}} \gtrsim\left\langle\xi_3^3\right\rangle^{\frac{1}{2}}$.
			\item \textbf{Case A1.2.2} If $\left(\xi_1, \xi_2, \xi_3\right)$ lies in \textbf{(b)}, that is, $\xi_1 \geq 0$,$ \xi_2 \leq 0$, $\xi_3 \leq 0$. Then, one has
	 $$H=\xi_2\left(3 \xi_1^2+3 \xi_1 \xi_2+2 \xi_2^2-1\right).$$
			Hence, it leads to 
			$
			|H| \gtrsim |\xi_2|^3 \sim|\xi_3|^3.
			$
			Moreover, we notice that, for $2-2b > \frac{2}{3}$ and $2b > \frac{2}{3}$,
			$
			\left\langle L_3\right\rangle^{2-2b}\left\langle L_1+L_2\right\rangle^{2 b} \gtrsim\left\langle H\right\rangle^{\frac{2}{3}} \sim\left\langle\xi_3\right\rangle^2 .
			$
			Therefore, one has,
			$$
			(\ref{3.6}) \lesssim \frac{\left|\xi_3\right|^2}{\left\langle\xi_3\right\rangle^2} \int_{\left|\xi_1\right| \leq 1} d \xi_1 \lesssim 1.
			$$
			
			\item \textbf{Case A1.2.3} If $\left(\xi_1, \xi_2, \xi_3\right)$ lies in \textbf{ (c)}, that is, $\xi_1 \leq 0$, $\xi_2 \geq 0$, $\xi_3 \leq 0$. Then, one has
			$$
			R_1\left(\xi_1\right) := L_1+L_2=\xi_1\left(2 \xi_1^2+3 \xi_1 \xi_3+3 \xi_3^2-1\right)-L_3,
			$$
			with $\left|R_1^{\prime}\left(\xi_1\right)\right| \gtrsim\left|\xi_1\right|^2+\left|\xi_1+\xi_3\right|^2 \gtrsim\left|\xi_2\right|^2
			$.
			Hence, we obtain that, for $\frac{1}{2} < b\leq b_0<1$,
			$$
			(\ref{3.6}) \lesssim \frac{\left|\xi_2\right|^2}{\left\langle L_3\right\rangle^{2-2b}} \int \frac{1}{\left|R_1^{\prime}\left(\xi_1\right)\right|} \frac{\left|R_1^{\prime}\left(\xi_1\right)\right|}{\left\langle R_1\left(\xi_1\right)\right\rangle^{2 b}} d \xi_1 \lesssim \int \frac{\left|R_1^{\prime}\left(\xi_1\right)\right|}{\left\langle R_1\left(\xi_1\right)\right\rangle^{2 b}} d \xi_1 \lesssim 1.
			$$
		\end{itemize}
	\end{itemize}
		\textbf{Case A2.} $\left|\xi_2\right| \leq 1$. This case can be addressed similarly as in  \textbf{Case A1.}\\
		\textbf{Case A3.} $\left|\xi_1\right|,\left|\xi_2\right|>1$ and $\left\langle L_1\right\rangle=\operatorname{MAX}:=\max \left\{\left\langle L_1\right\rangle,\left\langle L_2\right\rangle,\left\langle L_3\right\rangle\right\}$. This follows that
		$$
		\left\langle L_1\right\rangle^{b}\left\langle L_2\right\rangle^{b}\left\langle L_3\right\rangle^{1-b} \geq\left\langle L_1\right\rangle^{1- b}\left\langle L_2\right\rangle^b\left\langle L_3\right\rangle^b .
		$$
		Hence, to establish estimate \eqref{3.4}, it  suffices to bound the term,
		\begin{align}
			\sup _{\xi_1, \tau_1} \frac{\left\langle\xi_1\right\rangle^{2 p}}{\left\langle L_1\right\rangle^{2-2b}} \int \frac{\left|\xi_3\right|^2\left\langle\xi_2\right\rangle^{2 p}}{\left\langle\xi_3\right\rangle^{2 p}\left\langle L_2+L_3\right\rangle^{2 b}} d \xi_2,
			\label{3.7}
		\end{align}
	
	\begin{itemize}
		\item \textbf{Case A3.1.} If $\left(\xi_1, \xi_2, \xi_3\right)$ lies in \textbf{(a)}, that is, $\xi_1 \geq 0$, $\xi_2 \geq 0$, $\xi_3 \leq 0$, we observe that
		$$
		\begin{aligned}
		H  =  3 \xi_1 \xi_2\left(\xi_1+\xi_2\right),\quad 	R_2\left(\xi_2\right) & :=L_2+L_3
			 =3 \xi_1 \xi_2^2+3 \xi_1^2 \xi_2-L_1,
		\end{aligned}
		$$
		and $\left|R_2^{\prime}\left(\xi_2\right)\right|=|3 \xi_1\left(2 \xi_2+\xi_1\right)| \gtrsim\left|\xi_1 \xi_2\right|.$ Thus, it leads to 
		$
		\langle L_1\rangle \gtrsim \langle H \rangle\sim\left\langle \xi_1 \xi_2 \xi_3\right\rangle.
		$
		Then, for $ p < \frac{3}{4}$, one has, 
		\begin{align*}
			\left\langle\xi_1\right\rangle^{-1+2p}\left\langle\xi_2\right\rangle^{-1+2p}\left|\xi_3\right|^{2-2p} \lesssim\left|\xi_1 \xi_2 \xi_3\right|^{2-2 p} \lesssim\left\langle\xi_1 \xi_2 \xi_3\right\rangle^{2-2 p} \lesssim\left\langle L_1\right\rangle^{2-2 p},
		\end{align*}
		due to $\left|\xi_1\right|,\left|\xi_2\right|>1$ and    $-1+2p < 2-2p$.
		It follows that (\ref{3.7}) is bounded by
		\begin{align*}
			 \frac{\left\langle\xi_1\right\rangle^{2 p}}{\left\langle L_1\right\rangle^{2-2b}} \int \frac{\left|\xi_3\right|^2\left\langle\xi_2\right\rangle^{2 p}}{\left|\xi_3\right|^{2 p}\left|\xi_1 \xi_2\right|} \frac{\left|R_2^{\prime}\left(\xi_2\right)\right|}{\left\langle R_2\left(\xi_2\right)\right\rangle^{2 b}} d \xi_2 
			\lesssim  \frac{\left\langle L_1\right\rangle^{2-2 p}}{\left\langle L_1\right\rangle^{2-2b}} \lesssim 1
		\end{align*}
	with $2-2 p \leq 2-2b$.
		\item \textbf{Case A3.2.} If $\left(\xi_1, \xi_2, \xi_3\right)$ lies in \textbf{(b)} or \textbf{(c)}, the arguments are similar to the statements in \textbf{Case A3.1}.
	\end{itemize}
		 \textbf{Case A4.} $\left|\xi_1\right|,\left|\xi_2\right|>1$ and $\left\langle L_2\right\rangle=$ MAX. Due to the symmetric structure between $\xi_2$ and $\xi_3$ in (\ref{3.4}), the proof is identical to \textbf{Case A3.}\\
		 \textbf{Case A5.} $\left|\xi_1\right|,\left|\xi_2\right|>1$ and $\left\langle L_3\right\rangle=$ MAX.
		It is  suffices to bound the term,
		\begin{align}
			\sup _{\xi_3, \tau_3} \frac{\left|\xi_3\right|^2}{\left\langle\xi_3\right\rangle^{2 p}\left\langle L_3\right\rangle^{2-2b}} \int \frac{\left\langle\xi_1\right\rangle^{2 p}\left\langle\xi_2\right\rangle^{2 p}}{\left\langle L_1+L_2\right\rangle^{2 b}} d \xi_1,
			\label{3.8}
		\end{align}
	
	\begin{itemize}
		\item \textbf{Case A5.1.} If $\left(\xi_1, \xi_2, \xi_3\right)$ lies in \textbf{(a)}, that is, $\xi_1 \geq 0$, $\xi_2 \geq 0$, $\xi_3 \leq 0$. In this case, we have $|\xi_3| \geq |\xi_1|$, $|\xi_3| \geq |\xi_2|$,
		$$
		\begin{aligned}
			R_3\left(\xi_1\right)  :=L_1+L_2
			 =3 \xi_3 \xi_1^2+3 \xi_3^2 \xi_1-L_3,
		\end{aligned}
		$$
		and $R_3^{\prime}\left(\xi_1\right)=3 \xi_3\left(2 \xi_1+\xi_3\right)$. Moreover, we notice that
		$
		\left\langle L_3\right\rangle \gtrsim
		\left\langle H\right\rangle=\left\langle3 \xi_1 \xi_3\left(\xi_1+\xi_3\right)\right\rangle \gtrsim\left\langle\xi_1 \xi_2 \xi_3\right\rangle.
		$
		\begin{itemize}
			\item \textbf{Case A5.1.1}  If $\left|\xi_1\right| \lesssim  \left|\xi_3\right|$, then $\left|\xi_1\right| \lesssim\left|\xi_3\right| \sim\left|\xi_2\right|$ and $\left|R_3^{\prime}\left(\xi_1\right)\right| \gtrsim\left|\xi_3\right|^2$. 
			In addition, for $p>0$ and $|\xi_1|,|\xi_2|>1$, one has
			$$
			\left\langle\xi_1\right\rangle^{2 p}\left\langle\xi_2\right\rangle^{2 p}|\xi_3|^{-2p}
			\lesssim
			\left|\xi_1\right|^{2 p}\left|\xi_2\right|^{2 p}\left|\xi_3\right|^{-2 p} \lesssim\left|\xi_1\right|^{2 p} \lesssim\left|\xi_1 \xi_2 \xi_3\right|^{\frac{2}{3}p} \lesssim\left\langle L_3\right\rangle^{\frac{2}{3}p}.
			$$
			Thus,  (\ref{3.8}) is bounded by
             $$
				 \frac{\left|\xi_3\right|^2}{\left|\xi_3\right|^{2 p}\left\langle L_3\right\rangle^{2-2b}} \int \frac{\left\langle\xi_1\right\rangle^{2 p}\left\langle\xi_2\right\rangle^{2 p}}{\left|\xi_3\right|^2} \frac{\left|R_3^{\prime}\left(\xi_1\right)\right|}{\left\langle R_3\left(\xi_1\right)\right\rangle^{2 b}} d \xi_1
			\lesssim
			 \frac{\left\langle L_3\right\rangle^{\frac{2}{3}p}}{\left\langle L_3\right\rangle^{2-2b}},
			$$
			which is finite since $\frac{2}{3}p \leq 2-2b$.
			\item \textbf{Case A5.1.2} If $\left|\xi_1\right| \sim \left|\xi_3\right|$, then $\left|\xi_2\right| \lesssim\left|\xi_1\right| \sim\left|\xi_3\right|$. 
			We observe that, for $p>0$,
			$$
			\left\langle\xi_1\right\rangle^{2 p}\left\langle\xi_2\right\rangle^{2 p}\left\langle\xi_3\right\rangle^{\frac{3}{2}-2 p} \lesssim\left|\xi_1 \xi_2 \xi_3\right|^{\frac{1}{2}+\frac{2}{3}p} \lesssim\left\langle L_3\right\rangle^{\frac{1}{2}+\frac{2}{3}p},
			$$
			since $2 p \geq \frac{1}{2}+\frac{2}{3} p$. Additionally, according to  Lemma \ref{lem3},  (\ref{3.8}) is bounded by
			$$
			 \frac{\left|\xi_3\right|^{\frac{1}{2}}}{\left\langle L_3\right\rangle^{2-2b}} \int \frac{\left|\xi_3\right|^{\frac{3}{2}}\left\langle\xi_1\right\rangle^{2 p}\left\langle\xi_2\right\rangle^{2 p}}{\left|\xi_3\right|^{2 p}} \frac{1}{\left\langle R_3\left(\xi_1\right)\right\rangle^{2 b}} d \xi_1 \lesssim \frac{\left\langle L_3\right\rangle^{\frac{1}{2}+\frac{2}{3}p}}{\left\langle L_3\right\rangle^{2-2b}},
			$$
			which is finite for $\frac{1}{2}+\frac{2}{3}p \leq 2-2b$.
			\item \textbf{Case A5.1.3} If $\left|\xi_1\right| \gg \left|\xi_3\right|$, it implies that $|\xi_1|>|\xi_3|$. Therefore, the case does not exist for $|\xi_3|\geq|\xi_1|$.
		\end{itemize}
		\item \textbf{Case A5.2.}  If $\left(\xi_1, \xi_2, \xi_3\right)$ lies in \textbf{(b)} or \textbf{(c)}, these cases can be addressed same as \textbf{ Case A3.1}.
	\end{itemize}
	Next, we proceed to examine the case that $\left(\tau_1, \tau_2, \tau_3\right)$ lies in region $\mathbf{B}$. We  write,
	$$
	L_1=\tau_1-\left|\xi_1\right|^3+\frac{1}{2}\left|\xi_1\right|, \quad L_2=\tau_2+\left|\xi_2\right|^3-\frac{1}{2}\left|\xi_2\right|, \quad L_3=\tau_3-\left|\xi_3\right|^3+\frac{1}{2}\left|\xi_3\right| .
	$$
	\textbf{Case B1.}  If $\left|\xi_1\right| \leq 1$, similar to the statements in \textbf{Case A1}, it suffices to prove
		\begin{align*}
			\int_{A_3} \frac{\left|\xi_3\right| \prod_{j=1}^3\left|h_j\left(\xi_j, \tau_j\right)\right|}{\left\langle L_1\right\rangle^{b}\left\langle L_2\right\rangle^{b}\left\langle L_3\right\rangle^{1-b}} \lesssim \prod_{j=1}^3\left\|h_j\right\|_{L_{\xi,\tau}^2} .
		\end{align*}
	\begin{itemize}
		\item \textbf{Case B1.1} If $\left|\xi_2\right| \leq 2$, the proof   follows a similar approach as in \textbf{Case A1.1}.
		\item \textbf{Case B1.2} If $\left|\xi_2\right|>2$,   we have $\left|\xi_2\right| \sim\left|\xi_3\right|$ and it suffices to bound
		\begin{align}
			\sup _{\xi_3, \tau_3} \frac{\left|\xi_3\right|^2}{\left\langle L_3\right\rangle^{2-2b}} \int \frac{d \xi_1}{\left\langle L_1+L_2\right\rangle^{2 b}} .
			\label{3.9}
		\end{align}
		\item \textbf{Case B1.2.1} If $\left(\xi_1, \xi_2, \xi_3\right)$ lies in \textbf{(a)}, that is, $\xi_1 \geq 0$, $\xi_2 \geq 0$, $\xi_3 \leq 0$. Then, one has
		$$
		R_4\left(\xi_1\right) := L_1+L_2 =-\xi_1\left(2 \xi_1^2+3 \xi_1 \xi_3+3 \xi_3^2-1\right)-L_3,
		$$
		with $\left|R_4^{\prime}\left(\xi_1\right)\right| \gtrsim\left|\xi_3\right|^2$.
		Hence, we claim that, for $2-2b>0$,
		$$
		(\ref{3.9}) \lesssim   \frac{\left|\xi_3\right|^2}{\left\langle L_3\right\rangle^{2-2b}} \int \frac{1}{\left|R_4^{\prime}\left(\xi_1\right)\right|} \frac{\left|R_4^{\prime}\left(\xi_1\right)\right|}{\left\langle R_4\left(\xi_1\right)\right\rangle^{2 b}} d \xi_1 \lesssim 1 .
		$$
		\item \textbf{Case B1.2.2} If $\left(\xi_1, \xi_2, \xi_3\right)$ lies in \textbf{(b)}, that is, $\xi_1 \geq 0, \xi_2 \leq 0, \xi_3 \leq 0$. Then, one has
		$$|H|=\left|\xi_3\left(2 \xi_3^2+3 \xi_1 \xi_3+3 \xi_1^2-1\right)\right| 
		\sim\left|\xi_3\right|^3 .$$
		Hence, for $2-2b>\frac{2}{3}$ and $2b>\frac{2}{3}$, we have
		$
		\left\langle L_3\right\rangle^{2-2b}\left\langle L_1+L_2\right\rangle^{2 b} \gtrsim\left\langle H\right\rangle^{\frac{2}{3}}
		$,
	 the conclusion follows from a similar discussion as in \textbf{Case A1.2.2.}
		\item \textbf{Case B1.2.3} If $\left(\xi_1, \xi_2, \xi_3\right)$ lies in \textbf{(c)},
%
		 we can just repeat the proof in \textbf{Case A1.2.1.}
	\end{itemize}
		 \textbf{Case B2.} $\left|\xi_2\right| \leq 1$. The proof for this case is   similar to \textbf{Case B1}.
		\\
		 \textbf{Case B3.} $\left|\xi_1\right|,\left|\xi_2\right|>1$ and $\left\langle L_1\right\rangle=\operatorname{MAX}$.
		It suffices to bound
		\begin{align}
			\sup _{\xi_1, \tau_1} \frac{\left\langle\xi_1\right\rangle^{2 p}}{\left\langle L_1\right\rangle^{2-2b}} \int \frac{\left|\xi_3\right|^2\left\langle\xi_2\right\rangle^{2 p}}{\left\langle\xi_3\right\rangle^{2 p}\left\langle L_2+L_3\right\rangle^{2 b}} d \xi_2.
			\label{3.10}
		\end{align}
		Regardless of whether $\left(\xi_1, \xi_2, \xi_3\right)$ lies in \textbf{(a)}, \textbf{(b)} or \textbf{(c)}, the arguments for these   cases are   similar to  the one in \textbf{Case A3.1.}
		\\
		 \textbf{Case B4.}  $\left|\xi_1\right|,\left|\xi_2\right|>1$ and $\left\langle L_2\right\rangle=$ MAX. It yields that
		\begin{align*}
			\left\langle L_1\right\rangle^{b}\left\langle L_2\right\rangle^{b}\left\langle L_3\right\rangle^{1-b} \geq\left\langle L_2\right\rangle^{1-b}\left\langle L_1\right\rangle^b\left\langle L_3\right\rangle^b .
		\end{align*}
		Then, it suffices to bound
		\begin{align}
			\sup _{\xi_2, \tau_2} \frac{\left\langle\xi_2\right\rangle^{2 p}}{\left\langle L_2\right\rangle^{2-2b}} \int \frac{\left|\xi_3\right|^2\left\langle\xi_1\right\rangle^{2 p}}{\left\langle\xi_3\right\rangle^{2 p}\left\langle L_1+L_3\right\rangle^{2 b}} d \xi_1.
			\label{3.11}
		\end{align}
	\begin{itemize}
		\item \textbf{Case B4.1.}  If $\left(\xi_1, \xi_2, \xi_3\right)$ lies in \textbf{(a)} or \textbf{(b)}, this case can be addressed same as \textbf{ Case A3.1}.
		\item \textbf{Case B4.2.} If $\left(\xi_1, \xi_2, \xi_3\right)$ lies in \textbf{(c)}, that is, $\xi_1 \leq 0$, $\xi_2 \geq 0$, $\xi_3 \leq 0$.
		In this case, we have $|\xi_2|\geq|\xi_1|$, $|\xi_2|\geq|\xi_3|$,
		$$
		R_5\left(\xi_1\right):=L_1+L_3=-3 \xi_1 \xi_2\left(\xi_1+\xi_2\right)-L_2
		$$
		and
		$
		R_5^{\prime}\left(\xi_1\right)=-3 \xi_2\left(2 \xi_1+\xi_2\right)
		$.
		Moreover, one also has
		$
		\left\langle L_2\right\rangle \gtrsim
		\left\langle H\right\rangle=\left\langle3 \xi_1 \xi_2\left(\xi_1+\xi_2\right)\right\rangle \gtrsim\left\langle\xi_1 \xi_2 \xi_3\right\rangle.
		$
		\begin{itemize}
			\item \textbf{Case B4.3.1.} If $\left|2 \xi_1+\xi_2\right| \gtrsim \left|\xi_1\right|$, one obtains that $\left|R_5^{\prime}\left(\xi_2\right)\right| \gtrsim\left|\xi_1 \xi_2\right|$ and  the rest of the proof is similar to the one in \textbf{Case A3.1}.
			\item \textbf{Case B4.3.2.} If $\left|2 \xi_1+\xi_2\right| \ll \left|\xi_1\right|$, it implies that
			$\left|\xi_1\right| \sim\left|\xi_2\right| \sim\left|\xi_3\right|$. Then, we notice that
			$$
			\left\langle\xi_1\right\rangle^{2 p}\left\langle\xi_2\right\rangle^{-\frac{1}{2}+2p}\left|\xi_3\right|^{2-2p} \sim\left|\xi_1\right|^{\frac{3}{2}+2 p} \lesssim\left|\xi_1 \xi_2 \xi_3\right|^{\frac{1}{2}+\frac{2}{3}p} \lesssim\left\langle L_1\right\rangle^{\frac{1}{2}+\frac{2}{3}p},
			$$
			since $\left|\xi_1\right|,\left|\xi_2\right|>1$.
			Hence, according to Lemma \ref{lem3}, (\ref{3.11}) is bounded by
			\begin{align*}
				 \frac{\left\langle\xi_2\right\rangle^{\frac{1}{2}}}{\left\langle L_2\right\rangle^{2-2b}} \int \frac{\left|\xi_3\right|^2\left\langle\xi_2\right\rangle^{-\frac{1}{2}+2p}\left\langle\xi_1\right\rangle^{2 p}}{\left|\xi_3\right|^{2 p}} \frac{1}{\left\langle R_5\left(\xi_1\right)\right\rangle^{2 b}} d \xi_1 \lesssim
				 \frac{\left\langle L_2\right\rangle^{\frac{1}{2}+\frac{2}{3}p}}{\left\langle L_2\right\rangle^{2-2b}},
			\end{align*}
			which is finite since $\frac{1}{2}+\frac{2}{3}p \leq 2-2b$.
		\end{itemize}
	\end{itemize}
	\textbf{Case B5.} $\left|\xi_1\right|,\left|\xi_2\right|>1$ and $\left\langle L_3\right\rangle=$ MAX. It suffices to show the bound for
		\begin{align}
			\sup _{\xi_3, \tau_3} \frac{\left|\xi_3\right|^2}{\left\langle\xi_3\right\rangle^{2 p}\left\langle L_3\right\rangle^{2-2b}} \int \frac{\left\langle\xi_1\right\rangle^{2 p}\left\langle\xi_2\right\rangle^{2 p}}{\left\langle L_1+L_2\right\rangle^{2 b}} d \xi_1.
			\label{3.12}
		\end{align}
	\begin{itemize}
		\item \textbf{Case B5.1.}  If $\left(\xi_1, \xi_2, \xi_3\right)$ lies in \textbf{(a)}, we can follow a similar statements as in  \textbf{Case A3.1}.
		\item \textbf{Case B5.2.}  If $\left(\xi_1, \xi_2, \xi_3\right)$ lies in \textbf{(b)}, that is, $\xi_1 \geq 0$, $\xi_2 \leq 0$, $\xi_3 \leq 0$. In this case, we have $|\xi_1|\geq|\xi_2|$, $|\xi_1|\geq|\xi_3|$,
		$$
		R_6\left(\xi_1\right):=L_1+L_2=\xi_3\left(2 \xi_3^2+3 \xi_1 \xi_3+3 \xi_1^2-1\right) -L_3
		$$
		and
		$
		R_6^{\prime}\left(\xi_1\right)=3 \xi_3\left(2 \xi_1+\xi_3\right)
		$.
		Furthermore, we claim that
		$$
		\left\langle L_3\right\rangle \gtrsim
		\left\langle H\right\rangle=\left\langle \xi_3\left(2 \xi_3^2+3 \xi_1 \xi_3+3 \xi_1^2-1\right)\right\rangle \gtrsim\left\langle\xi_1 \xi_2 \xi_3\right\rangle.
		$$
		To bound (\ref{3.12}), we can just repeat the proof in \textbf{Case A5.1}.
		\item \textbf{Case B5.3.}  If $\left(\xi_1, \xi_2, \xi_3\right)$ lies in \textbf{(c)}, that is, $\xi_1 \leq 0$, $\xi_2 \geq 0$, $\xi_3 \leq 0$. In this case, we have $|\xi_2|\geq|\xi_1|$, $|\xi_2|\geq|\xi_3|$,
	$$
		R_7\left(\xi_1\right):=L_1+L_2=-3 \xi_1 \xi_3\left(\xi_1+\xi_3\right)-L_3
		$$
		and
		$
		R_7^{\prime}\left(\xi_1\right)=-3 \xi_3\left(2 \xi_1+\xi_3\right)
		$.
		Additionally, we point out that
		$$
		\left\langle L_3\right\rangle \gtrsim
		\left\langle H\right\rangle=\left\langle 3 \xi_1 \xi_3\left(\xi_1+\xi_3\right)\right\rangle \gtrsim\left\langle\xi_1 \xi_2 \xi_3\right\rangle.
		$$
		\begin{itemize}
			\item \textbf{Case B5.3.1.} If $\left|\xi_1\right| \ll\left|\xi_3\right|$, then $\left|\xi_1\right| \lesssim\left|\xi_3\right| \sim\left|\xi_2\right|$ and $\left|R_7^{\prime}\left(\xi_1\right)\right| \gtrsim\left|\xi_3\right|^2$.
			We notice that, for $p>0$,
			$$
			\left\langle\xi_1\right\rangle^{2 p}\left\langle\xi_2\right\rangle^{2 p}\left|\xi_3\right|^{-2 p} \lesssim\left|\xi_1\right|^{2 p} \lesssim\left|\xi_1 \xi_2 \xi_3\right|^{\frac{2}{3}p} \lesssim\left\langle L_3\right\rangle^{\frac{2}{3}p},
			$$
			since $\left|\xi_1\right|,\left|\xi_2\right|>1$.
			 Thus, (\ref{3.12}) is bounded by
			\begin{align*}
				 \frac{\left|\xi_3\right|^2}{\left|\xi_3\right|^{2 p}\left\langle L_3\right\rangle^{2-2b}} \int \frac{\left\langle\xi_1\right\rangle^{2 p}\left\langle\xi_2\right\rangle^{2 p}}{\left|\xi_3\right|^2}
				\frac{\left|R_7^{\prime}\left(\xi_1\right)\right|}{\left\langle R_7\left(\xi_1\right)\right\rangle^{2 b}} d \xi_1 \lesssim
				 \frac{\left\langle L_3\right\rangle^{\frac{2}{3}p}}{\left\langle L_3\right\rangle^{2-2b}},
			\end{align*}
			which is finite since $\frac{2}{3}p \leq 2-2b$.
			\item \textbf{Case B5.3.2.}  If $ \left|\xi_1\right| \sim \left|\xi_3\right|$, then $\left|\xi_2\right| \lesssim\left|\xi_1\right| \sim\left|\xi_3\right|$.  
			We notice that, for $p>0$,
			$$
			\left\langle\xi_1\right\rangle^{2 p}\left\langle\xi_2\right\rangle^{2 p}\left|\xi_3\right|^{\frac{3}{2}-2 p} \lesssim\left|\xi_1 \xi_2 \xi_3\right|^{\frac{1}{2}+\frac{2}{3}p} \lesssim\left\langle L_3\right\rangle^{\frac{1}{2}+\frac{2}{3}p},
			$$
			since $2 p \geq \frac{1}{2}+\frac{2}{3} p$ and $\left|\xi_1\right|,\left|\xi_2\right|>1$. Thus, according to Lemma \ref{lem3}, (\ref{3.12}) is bounded by
			\begin{align*}
				 \frac{\left|\xi_3\right|^{\frac{1}{2}}}{\left\langle L_3\right\rangle^{2-2b}} \int \frac{\left|\xi_3\right|^{\frac{3}{2}}\left\langle\xi_1\right\rangle^{2 p}\left\langle\xi_2\right\rangle^{2 p}}{\left|\xi_3\right|^{2p}}
				\frac{1}{\left\langle R_7\left(\xi_1\right)\right\rangle^{2 b}} d \xi_1 \lesssim
				 \frac{\left\langle L_3\right\rangle^{\frac{1}{2}+\frac{2}{3}p}}{\left\langle L_3\right\rangle^{2-2b}},
			\end{align*}
			which is finite for $\frac{1}{2}+\frac{2}{3}p \leq 2-2b$.
			\item \textbf{Case B5.3.3.} If $\left|\xi_1\right| \gg\left|\xi_3\right|$, then $\left|\xi_3\right| \lesssim\left|\xi_1\right| \sim\left|\xi_2\right|$ and $\left|R_7^{\prime}\left(\xi_1\right)\right| \gtrsim\left|\xi_1 \xi_3\right|$. 
			We notice that, for $\frac{9}{16}<p<\frac{3}{4}$, 
			one has
			$$
			\left\langle\xi_1\right\rangle^{-1+2p}\left\langle\xi_2\right\rangle^{2 p}\left\langle\xi_3\right\rangle^{-2 p} \left|\xi_3\right|\lesssim\left|\xi_1 \xi_2 \xi_3\right|^{-\frac{1}{2}+2p} \lesssim\left\langle L_3\right\rangle^{-\frac{1}{2}+2p}.
			$$
			Thus,  (\ref{3.12}) is bounded by
			$$
			 \frac{\left|\xi_3\right|^2}{\left\langle\xi_3\right\rangle^{2 p}\left\langle L_3\right\rangle^{2-2b}} \int \frac{\left\langle\xi_1\right\rangle^{2 p}\left\langle\xi_2\right\rangle^{2 p}}{\left|\xi_1 \xi_3\right|} \frac{\left|R_7^{\prime}\left(\xi_1\right)\right|}{\left\langle R_7\left(\xi_1\right)\right\rangle^{2 b}} d \xi_1 \lesssim
			 \frac{\left\langle L_3\right\rangle^{-\frac{1}{2}+2p}}{\left\langle L_3\right\rangle^{2-2b}},
			$$
			which is finite if $-\frac{1}{2}+2p \leq 2-2b$.
		\end{itemize}
	\end{itemize}
	The proof is now complete. \\
	\subsection{Proof of Estimate \eqref{bi13}}
 Similar to the setting in the proof of \textbf{Proposition \ref{bi0}},
	based on duality and Plancherel theorem, \eqref{bi13} is equivalent to
	\begin{align}
		\int_{A_3} \left|\frac{|\xi_3|\xi_2^2 h_1(\xi_1,\tau_1)h_2(\xi_2,\tau_2)h_3\left(\xi_3, \tau_3\right)}{\left\langle\xi_1\right\rangle^s\left\langle\xi_2\right\rangle^s\left\langle\xi_3\right\rangle^{2-s}\left\langle L_1\right\rangle^{b}\left\langle L_2\right\rangle^{b}\left\langle L_3\right\rangle^{1-b}} \right|\lesssim \prod_{j=1}^3\left\|h_j\right\|_{L_{\xi, \tau}^2}. \label{3.13}
	\end{align}
	Since we only need to consider $s$ close to $\frac{1}{2}$, without loss of generality, we assume $\frac{1}{2}<s<\frac{5}{8}$ and define $b_0 = s$ to help to justify the proof for $\frac{1}{2}<b \leq b_0$. It again suffices to establish the estimate for  Cases \textbf{(a)}, \textbf{(b)} and \textbf{(c)} lying in sets \textbf{A} and \textbf{B}. 
	
	\textbf{Case A1.} If $\left|\xi_1\right| \leq 1$, since $\left\langle\xi_1\right\rangle\left\langle\xi_2\right\rangle /\left\langle\xi_3\right\rangle \lesssim 1$, estimate (\ref{3.13}) can be simplified to
		\begin{align}
			\int_{A_3(\mathbb{R}) } \frac{\left|\xi_2\right|^2 \prod_{j=1}^3\left|h_j\left(\xi_j, \tau_j\right)\right|}{\left\langle\xi_3\right\rangle \left\langle L_1\right\rangle^{b}\left\langle L_2\right\rangle^{b}\left\langle L_3\right\rangle^{1-b}} \lesssim \prod_{j=1}^3\left\|h_j\right\|_{L_{\xi, \tau}^2}.
			\label{3.14}
		\end{align}
	\begin{itemize}
		\item \textbf{Case A1.1.}  If $\left|\xi_2\right| \leq 2$, one has $\left|\xi_3\right| \leq 3$. The proof is exactly same as \textbf{Case A1.1} in the proof of {\bf Proposition \ref{bi0}}.
		
		\item \textbf{Case A1.2.}  If $\left|\xi_2\right|>2$,
		we notice that $|\xi_1| \leq 1$ and $|\xi_2| > 2$ will lead to $|\xi_2| \sim |\xi_3|$. Thus, one has
		\begin{align*}
			\frac{1}{\langle\xi_3\rangle^2\left\langle L_3\right\rangle^{2-2b}} \int \frac{\left|\xi_2\right|^4d \xi_1}{\left\langle L_1+L_2\right\rangle^{2 b}}
			\sim \frac{\left|\xi_3\right|^2}{\left\langle L_3\right\rangle^{2-2b}} \int \frac{d \xi_1}{\left\langle L_1+L_2\right\rangle^{2 b}}.
		\end{align*}
		Hence, the proof of (\ref{3.14}) is equivalent to \textbf{Case A1.2} in the proof of {\bf Proposition \ref{bi0}} with $b \leq b_0<\frac{5}{8}$.

	\end{itemize}
	\textbf{Case A2.}  $\left|\xi_2\right| \leq 1$. 
		 This case can be addressed similar to  \textbf{Case A1.1}.\\
		 \textbf{Case  A3.} $\left|\xi_1\right|>1, \left|\xi_2\right|>1, \left|\xi_3\right| \leq 1$. One has $|\xi_1| \sim |\xi_2|$. For $\frac{1}{2} <b <1$, it suffices to bound the term
		 \begin{align*}
		 	\sup _{\xi_1, \tau_1} \frac{1}{\left\langle L_1\right\rangle^{2-2b}\left\langle\xi_1\right\rangle^{2s}} \int \frac{\left|\xi_2\right|^4}{\left\langle\xi_2\right\rangle^{2s}\left\langle L_2+L_3\right\rangle^{2-2 b}} d \xi_3.
		 \end{align*}
		 
		 \begin{itemize}
		 	\item \textbf{Case  A3.1} If $\left(\xi_1, \xi_2, \xi_3\right)$ lies in \textbf{(a)}, that is, $\xi_1 \geq 0, \xi_2 \geq 0, \xi_3 \leq 0$. It implies that $|\xi_3| = |\xi_1|+|\xi_2|$. Thus, this case does not exist. 
		 	
		 	\item \textbf{Case  A3.2} If $\left(\xi_1, \xi_2, \xi_3\right)$ lies in \textbf{(b)}, that is, $\xi_1 \geq 0, \xi_2 \leq 0, \xi_3 \leq 0$. Then, we have
		 	$$
		 	\left|H\right| = |\xi_2\left(3 \xi_1^2+3 \xi_1 \xi_2+2 \xi_2^2-1\right)| \gtrsim
		 	|\xi_1|^3.
		 	$$
		 	In addition, we notice that, for $\frac{1}{2}<b<1$,
		 	$
		 	\left\langle L_1\right\rangle^{2-2b}\left\langle L_2+L_3\right\rangle^{2-2 b} \gtrsim\left\langle H\right\rangle^{2-2b} \sim\left\langle\xi_1\right\rangle^{6-6b} .
		 	$
		 	Therefore, for $2-6b+4s>0$, one has
		 	$$
		 	\frac{1}{\left\langle L_1\right\rangle^{2-2b}\left\langle\xi_1\right\rangle^{2s}} \int \frac{\left|\xi_2\right|^4}{\left\langle\xi_2\right\rangle^{2s}\left\langle L_2+L_3\right\rangle^{2-2 b}} d \xi_3 \lesssim \frac{\left|\xi_1\right|^{4-4s}}{\left\langle\xi_1\right\rangle^{6-6b}} \int_{\left|\xi_3\right| \leq 1} d \xi_3 \lesssim 1.
		 	$$
		 	
		 	\item \textbf{Case  A3.3} If $\left(\xi_1, \xi_2, \xi_3\right)$ lies in \textbf{(c)}, that is, $\xi_1 \leq 0, \xi_2 \geq 0, \xi_3 \leq 0$. Then, one has
		 	$$
		 	\left|H\right| = |\xi_1\left(2 \xi_1^2+3 \xi_1 \xi_3+3 \xi_3^2-1\right)| \gtrsim
		 	|\xi_1|^3.
		 	$$
		 	Hence, the remainder proof is similar to the line in \textbf{Case  A3.2}.
		 	
		 \end{itemize}
	\textbf{Case A4.} $\left|\xi_1\right|,\left|\xi_2\right|,\left|\xi_3\right|>1$ and $\left\langle L_1\right\rangle=\operatorname{MAX}:=\max \left\{\left\langle L_1\right\rangle,\left\langle L_2\right\rangle,\left\langle L_3\right\rangle\right\}$.
		This follows that
		$$
		\left\langle L_1\right\rangle^{b}\left\langle L_2\right\rangle^{b}\left\langle L_3\right\rangle^{1-b} \geq\left\langle L_1\right\rangle^{1-b}\left\langle L_2\right\rangle^b\left\langle L_3\right\rangle^b .
		$$
		It then suffices to bound the term
		\begin{align}
			\sup _{\xi_1, \tau_1} \frac{1}{\left\langle L_1\right\rangle^{2-2b}\left\langle\xi_1\right\rangle^{2s}} \int \frac{\left|\xi_2\right|^4}{\left\langle\xi_2\right\rangle^{2s}\left\langle\xi_3\right\rangle^{2-2s}\left\langle L_2+L_3\right\rangle^{2 b}} d \xi_2,
			\label{3.15}
		\end{align}
	\begin{itemize}
		\item \textbf{Case A4.1.} If $\left(\xi_1, \xi_2, \xi_3\right)$ lies in \textbf{(a)}, that is, $\xi_1 \geq 0$, $\xi_2 \geq 0$, $\xi_3 \leq 0$. In this case, we observe that $|\xi_3|>|\xi_1|$, $|\xi_3|>|\xi_2|$,
		$$
		\begin{aligned}
			P_1\left(\xi_2\right) & :=L_2+L_3=3 \xi_1 \xi_2^2+3 \xi_1^2 \xi_2-L_1,
		\end{aligned}
		$$
		and $\left|P_1^{\prime}\left(\xi_2\right)\right| \gtrsim\left|\xi_1 \xi_2\right|$. Moreover, one has
		$
		\left\langle L_1\right\rangle \gtrsim\left\langle H\right\rangle = \langle3 \xi_1 \xi_2\left(\xi_1+\xi_2\right)\rangle \gtrsim\left\langle\xi_1 \xi_2 \xi_3\right\rangle.
		$
		 For $\frac{1}{2}<s<\frac{5}{8}$, it yields $-1+2s<2-2s$ and we have
			\begin{align*}
				\left\langle\xi_1\right\rangle^{-1-2s}\left\langle\xi_2\right\rangle^{3-2s}\left|\xi_3\right|^{-2+2s} \lesssim
				\left\langle\xi_1\right\rangle^{-1-2s}\left\langle\xi_2\right\rangle^{2-2s}\left|\xi_3\right|^{-1+2s} \lesssim\left|\xi_1 \xi_2 \xi_3\right|^{2-2s} \lesssim\left\langle\xi_1 \xi_2 \xi_3\right\rangle^{2-2s} \lesssim\left\langle L_1\right\rangle^{2-2s}.
			\end{align*}
			It follows that  (\ref{3.15})   is bounded by
			\begin{align*}
				 \frac{1}{\left\langle L_1\right\rangle^{2-2b}\left\langle\xi_1\right\rangle^{2s}} \int \frac{\left|\xi_2\right|^4}{\left\langle\xi_2\right\rangle^{2s}\left|\xi_3\right|^{2-2s}\left|\xi_1 \xi_2\right|} \frac{\left|P_1^{\prime}\left(\xi_2\right)\right|}{\left\langle P_1\left(\xi_2\right)\right\rangle^{2 b}} d \xi_2 \lesssim  \frac{\left\langle L_1\right\rangle^{2-2s}}{\left\langle L_1\right\rangle^{2-2b}},
			\end{align*}	
			which is finite owing to $2-2s \leq 2-2b$.
					
		\item \textbf{Case A4.2.}  If $\left(\xi_1, \xi_2, \xi_3\right)$ lies in \textbf{(b)}, that is, $\xi_1 \geq 0$, $\xi_2 \leq 0$, $\xi_3 \leq 0$. In this case, one has $|\xi_1|>|\xi_2|$, $|\xi_1|>|\xi_3|$,
		$$
		P_2\left(\xi_2\right):=L_2+L_3=2 \xi_2^3+3 \xi_1 \xi_2^2+\left(3 \xi_1^2-1\right) \xi_2-L_1,
		$$
		and
		$\left|P_2^{\prime}\left(\xi_2\right)\right| \gtrsim \xi_1^2 \gtrsim\left|\xi_1 \xi_2\right|$. Moreover, we also have
		$
		\left\langle L_1\right\rangle \gtrsim\left\langle H\right\rangle \gtrsim\left\langle\xi_1 \xi_2 \xi_3\right\rangle.
		$ Therefore, the rest of the discussion will be similar to  \textbf{Case A4.1.}
		\item \textbf{Case A4.3.} If $\left(\xi_1, \xi_2, \xi_3\right)$ lies in \textbf{(c)}, that is, $\xi_1 \leq 0, \xi_2 \geq 0, \xi_3 \leq 0$.
		In this case, we have $|\xi_2|>|\xi_1|$, $|\xi_2|>|\xi_3|$,
		$$
		P_3\left(\xi_2\right):=L_2+L_3=3 \xi_1 \xi_2^2+3 \xi_1^2 \xi_2-\tau_1+\xi_1^3-\frac{1}{2} \xi_1, \quad \left|P_3^{\prime}\left(\xi_2\right)\right| \gtrsim\left|\xi_1 \xi_2\right|,
		$$
		and
		$
		\left\langle L_1\right\rangle \gtrsim\left\langle H\right\rangle \gtrsim\left\langle\xi_1 \xi_2 \xi_3\right\rangle. 
		$ Then, we can follow the idea in \textbf{Case A4.1.}
	\end{itemize}
	 \textbf{Case A5.}  $\left|\xi_1\right|,\left|\xi_2\right|,\left|\xi_3\right|>1$ and $\left\langle L_2\right\rangle=\operatorname{MAX}$.
	This case can be addressed similar to  \textbf{Case A4}.
\\
	\textbf{Case A6.} $\left|\xi_1\right|,\left|\xi_2\right|,\left|\xi_3\right|>1$ and $\left\langle L_3\right\rangle=$ MAX.
		To obtain (\ref{3.13}), it is then sufficient to bound the term
		\begin{align}
			\sup _{\xi_3, \tau_3} \frac{1}{\left\langle\xi_3\right\rangle^{2-2s}\left\langle L_3\right\rangle^{2-2b}} \int \frac{\left|\xi_2\right|^4}{\left\langle\xi_1\right\rangle^{2s}\left\langle\xi_2\right\rangle^{2s}\left\langle L_1+L_2\right\rangle^{2 b}} d \xi_1.
			\label{3.17}
		\end{align}
	\begin{itemize}
		\item \textbf{Case A6.1.} If $\left(\xi_1, \xi_2, \xi_3\right)$ lies in \textbf{(a)}, that is, $\xi_1 \geq 0$, $\xi_2 \geq 0$, $\xi_3 \leq 0$. In this case, we have $|\xi_3| \geq |\xi_1|$, $|\xi_3| \geq |\xi_2|$,
		$$
		\begin{aligned}
			P_4\left(\xi_1\right)  :=L_1+L_2
			 =3 \xi_3 \xi_1^2+3 \xi_3^2 \xi_1-L_3,
		\end{aligned}
		$$
		and
		$
		P_4^{\prime}\left(\xi_1\right)=3 \xi_3\left(2 \xi_1+\xi_3\right).
		$
		Moreover, we notice that
		$
		\left\langle L_3\right\rangle \gtrsim\left\langle H\right\rangle = \langle3 \xi_1 \xi_3\left(\xi_1+\xi_3\right)\rangle \gtrsim\left\langle\xi_1 \xi_2 \xi_3\right\rangle.
		$
		\begin{itemize}
			\item \textbf{Case A6.1.1}  If $\left|\xi_1\right| \ll\left|\xi_3\right|$, then $\left|\xi_1\right| \lesssim\left|\xi_3\right| \sim\left|\xi_2\right|$ and $|P_4^{\prime}\left(\xi_1\right)| \gtrsim |\xi_3|^2$. 
			Thus, for $\frac{1}{2}<s<\frac{5}{8}$, one has,
			\begin{align*}
				\left\langle\xi_1\right\rangle^{-2s}\left\langle\xi_2\right\rangle^{4-2s}\left|\xi_3\right|^{-4+2s} \lesssim
				\left\langle\xi_1\right\rangle^{-2s} \lesssim\left|\xi_1 \xi_2 \xi_3\right|^{2-2s} \lesssim\left\langle\xi_1 \xi_2 \xi_3\right\rangle^{2-2s} \lesssim\left\langle L_1\right\rangle^{2-2s} ,
			\end{align*}
			Therefore,  (\ref{3.17}) is bounded by
			$$
			 \frac{1}{\left\langle L_3\right\rangle^{2-2b}} \int \frac{\left\langle\xi_2\right\rangle^{4-2s}}{\left\langle\xi_1\right\rangle^{2s}\left|\xi_3\right|^{4-2s}} \frac{|P_4^{\prime}\left(\xi_1\right)|}{\left\langle P_4\left(\xi_1\right)\right\rangle^{2 b}} d \xi_1
			\lesssim
			 \frac{\left\langle L_3\right\rangle^{2-2s}}{\left\langle L_3\right\rangle^{2-2b}},
			$$
			which is finite owing to $2-2s \leq 2-2b$.
			
			\item \textbf{Case A6.1.2} If $ \left|\xi_1\right| \sim\left|\xi_3\right|$, then $\left|\xi_2\right| \lesssim\left|\xi_1\right| \sim\left|\xi_3\right|$. For $\frac{1}{2}<s<\frac{5}{8}$, it yields $-\frac{1}{2}+2s<2-2s$ and we have
			$$
			\left\langle\xi_1\right\rangle^{-2s}\left\langle\xi_2\right\rangle^{4-2s}\left|\xi_3\right|^{-\frac{5}{2}+2s} \lesssim \left\langle\xi_1\right\rangle^{1-2s}\left\langle\xi_2\right\rangle^{1-2s}\left|\xi_3\right|^{-\frac{1}{2}+2s} \lesssim \left|\xi_1 \xi_2 \xi_3\right|^{2-2s} \lesssim\left\langle L_3\right\rangle^{2-2s}.
			$$
			Therefore, by using Lemma \ref{lem3}, (\ref{3.17}) is bounded by
			\begin{align*}
				 \frac{\left|\xi_3\right|^{\frac{1}{2}}}{\left\langle L_3\right\rangle^{2-2b}} \int \frac{\left\langle\xi_2\right\rangle^{4-2s}}{\left\langle\xi_1\right\rangle^{2s}\left|\xi_3\right|^{\frac{5}{2}-2s}} \frac{1}{\left\langle P_4\left(\xi_1\right)\right\rangle^{2 b}} d \xi_1
				\lesssim
				 \frac{\left\langle L_3\right\rangle^{2-2s}}{\left\langle L_3\right\rangle^{2-2b}},
			\end{align*}
			which is finite for $2-2s\leq 2-2b$.
			\item \textbf{Case A6.1.3} If $\left|\xi_1\right| \gg \left|\xi_3\right|$, it implies that $|\xi_1|>|\xi_3|$. Therefore, the case does not exist for $|\xi_3|\geq|\xi_1|$.
		\end{itemize}
		\item \textbf{Case A6.2.}  If $\left(\xi_1, \xi_2, \xi_3\right)$ lies in \textbf{(b)}, it is similar to the\textbf{ Case A4.2}. 
		\item \textbf{Case A6.3.} If $\left(\xi_1, \xi_2, \xi_3\right)$ lies in \textbf{(c)}, it is similar to the\textbf{ Case A4.3}.
		\end{itemize}
		Next, we proceed to examine the case that $\left(\tau_1, \tau_2, \tau_3\right)$ lies in region $\mathbf{B}$. 
		
	\noindent\textbf{Case B1.} $\left|\xi_1\right| \leq 1$.   
		The argument for this case will be similar to the one in \textbf{Case A1}.
		
	\noindent\textbf{Case B2.} $\left|\xi_2\right| \leq 1$. 
		In this case, we have $\left\langle\xi_1\right\rangle\left\langle\xi_2\right\rangle /\left\langle\xi_3\right\rangle \lesssim 1$ and $|\xi_2|^2/\langle\xi_3\rangle \lesssim 1.$
		This case can be addressed in the same way as \textbf{Case A2}.
		\\
		\textbf{Case B3.} $\left|\xi_1\right|>1, \left|\xi_2\right|>1, \left|\xi_3\right| \leq 1$. Thus, one has $|\xi_1| \sim |\xi_2|$. For $\frac{1}{2}<b<1$, it suffices to bound the term
		\begin{align*}
			\sup _{\xi_1, \tau_1} \frac{1}{\left\langle L_1\right\rangle^{2-2b}\left\langle\xi_1\right\rangle^{2s}} \int \frac{\left|\xi_2\right|^4}{\left\langle\xi_2\right\rangle^{2s}\left\langle L_2+L_3\right\rangle^{2-2 b}} d \xi_3.
		\end{align*}
		
		\begin{itemize}
			\item \textbf{Case B3.1} If $\left(\xi_1, \xi_2, \xi_3\right)$ lies in \textbf{(a)}, that is, $\xi_1 \geq 0, \xi_2 \geq 0, \xi_3 \leq 0$. It implies that $|\xi_3| = |\xi_1|+|\xi_2|$. Thus, this case does not exist. 
			
			\item \textbf{Case B3.2} If $\left(\xi_1, \xi_2, \xi_3\right)$ lies in \textbf{(b)}, that is, $\xi_1 \geq 0, \xi_2 \leq 0, \xi_3 \leq 0$. Then, we have
			$$
			\left|H\right| = |\left(\xi_1+\xi_2\right)\left(2 \xi_2^2+\xi_1 \xi_2+2 \xi_1^2-1\right)| \gtrsim
			|\xi_1\xi_2\xi_3|.
			$$
			In addition, we notice that, for $\frac{1}{2}<b<1$,
			$$
			\left\langle L_1\right\rangle^{2-2b}\left\langle L_2+L_3\right\rangle^{2-2 b} \gtrsim\left\langle H\right\rangle^{2-2b} \gtrsim\left\langle\xi_1\xi_2\xi_3\right\rangle^{2-2b} .
			$$
			Therefore, for $b-s<0$, one has,
			$$
			\frac{1}{\left\langle L_1\right\rangle^{2-2b}\left\langle\xi_1\right\rangle^{2s}} \int \frac{\left|\xi_2\right|^4}{\left\langle\xi_2\right\rangle^{2s}\left\langle L_2+L_3\right\rangle^{2-2 b}} d \xi_3 \lesssim \int_{\left|\xi_3\right| \leq 1} \frac{\left|\xi_1\right|^{4-4s}}{\left\langle\xi_1\xi_2\xi_3\right\rangle^{2-2b}} d \xi_3 \lesssim |\xi_1|^{4(b-s)} \lesssim 1.
			$$
			\item \textbf{Case B3.3} If $\left(\xi_1, \xi_2, \xi_3\right)$ lies in \textbf{ (c)}, that is, $\xi_1 \leq 0, \xi_2 \geq 0, \xi_3 \leq 0$. Then, one has
			$
			\left|H\right| = |3\xi_1\xi_2(\xi_1+\xi_2)| \sim
			|\xi_1\xi_2\xi_3|.
			$
			Hence, the remainder proof is similar to the line in \textbf{Case B3.2}.
			
		\end{itemize}
		\textbf{Case B4.} $\left|\xi_1\right|,\left|\xi_2\right|,\left|\xi_3\right|>1$ and $\left\langle L_1\right\rangle=$ MAX. 
	 This case can be addressed same as \textbf{Case A4.}\\
	\textbf{Case B5.} $\left|\xi_1\right|,\left|\xi_2\right|,\left|\xi_3\right|>1$ and $\left\langle L_2\right\rangle=$ MAX. It suffices to bound
		\begin{align}
			\sup _{\xi_2, \tau_2} \frac{1}{\left\langle L_2\right\rangle^{2-2b}\left\langle\xi_2\right\rangle^{2s}} \int \frac{\left|\xi_2\right|^4}{\left\langle\xi_1\right\rangle^{2s}\left\langle\xi_3\right\rangle^{2-2s}\left\langle L_1+L_3\right\rangle^{2 b}} d \xi_1
			\label{3.19}
		\end{align}
	\begin{itemize}
		\item \textbf{Case B5.1.}  If $\left(\xi_1, \xi_2, \xi_3\right)$ lies in \textbf{(a)}, it is similar to \textbf{Case A4.1}.
		\item \textbf{Case B5.2.}  If $\left(\xi_1, \xi_2, \xi_3\right)$ lies in \textbf{(b)}, it is similar to \textbf{Case A4.2}.
		
		\item \textbf{Case B5.3.} If $\left(\xi_1, \xi_2, \xi_3\right)$ lies in \textbf{(c)}, that is, $\xi_1 \leq 0$, $\xi_2 \geq 0$, $\xi_3 \leq 0$.
		In this case, we have $|\xi_2|\geq|\xi_1|$, $|\xi_2|\geq|\xi_3|$,
		$$
		P_5\left(\xi_1\right):=L_1+L_3=-3 \xi_1 \xi_2\left(\xi_1+\xi_2\right)-L_2,
		$$
		and
		$
		P_5^{\prime}\left(\xi_1\right)=-3 \xi_2\left(2 \xi_1+\xi_2\right)
		$.
		Moreover, one also has
		$
		\left\langle L_2\right\rangle \gtrsim\left\langle H\right\rangle = \langle3 \xi_1 \xi_2\left(\xi_1+\xi_2\right)\rangle \gtrsim\left\langle\xi_1 \xi_2 \xi_3\right\rangle.
		$
		\begin{itemize}
			\item \textbf{Case B5.3.1.} If $\left|2 \xi_1+\xi_2\right| \gg\left|\xi_1\right|$, then $\left|P_5^{\prime}\left(\xi_2\right)\right| \gtrsim\left|\xi_1 \xi_2\right|$. Thus, the rest of the proof is similar to \textbf{Case A4.3}.
			\item \textbf{Case B5.3.2.} If $\left|2 \xi_1+\xi_2\right| \lesssim\left|\xi_1\right|$, it implies that
			$\left|\xi_1\right| \sim\left|\xi_2\right| \sim\left|\xi_3\right|$. For $\frac{1}{2}<s<\frac{5}{8}$, it yields $-\frac{1}{2}+2s<2-2s$ and one has
			$$
			\left\langle\xi_1\right\rangle^{-2s}\left\langle\xi_2\right\rangle^{\frac{7}{2}-2s}\left\langle\xi_3\right\rangle^{-2+2s} \lesssim \left\langle\xi_1\right\rangle^{1-2s}\left\langle\xi_2\right\rangle^{1-2s}\left\langle\xi_3\right\rangle^{-\frac{1}{2}+2s} \lesssim \left|\xi_1 \xi_2 \xi_3\right|^{2-2s} \lesssim\left\langle L_3\right\rangle^{2-2s},
			$$
			since $|\xi_1|,|\xi_2|,|\xi_3|>1$.
				Therefore, by using Lemma \ref{lem3}, (\ref{3.19}) is bounded by
			\begin{align*}
				 \frac{\left|\xi_2\right|^{\frac{1}{2}}}{\left\langle L_2\right\rangle^{2-2b}} \int \frac{\left\langle\xi_2\right\rangle^{\frac{7}{2}-2s}}{\left\langle\xi_1\right\rangle^{2s}\left\langle\xi_3\right\rangle^{2-2s}} \frac{1}{\left\langle P_5\left(\xi_1\right)\right\rangle^{2 b}} d \xi_1
				\lesssim
				 \frac{\left\langle L_2\right\rangle^{2-2s}}{\left\langle L_2\right\rangle^{2-2b}},
			\end{align*}
			which is finite for $2-2s\leq 2-2b$.
		\end{itemize}
	\end{itemize}
	\textbf{Case B6.} $\left|\xi_1\right|,\left|\xi_2\right|,\left|\xi_3\right|>1$ and $\left\langle L_3\right\rangle=$ MAX. It suffices to show the bound for
		\begin{align}
			\sup _{\xi_3, \tau_3} \frac{1}{\left\langle\xi_3\right\rangle^{2-2s}\left\langle L_3\right\rangle^{2-2b}} \int \frac{\left|\xi_2\right|^4}{\left\langle\xi_1\right\rangle^{2s}\left\langle\xi_2\right\rangle^{2s}\left\langle L_1+L_2\right\rangle^{2 b}} d \xi_1,
			\label{3.20}
		\end{align}
	\begin{itemize}
		\item \textbf{Case B6.1.}  If $\left(\xi_1, \xi_2, \xi_3\right)$ lies in \textbf{(a)}, it is similar to \textbf{Case A4.1}.
		\item \textbf{Case B6.2.}  If $\left(\xi_1, \xi_2, \xi_3\right)$ lies in \textbf{(b)}, that is, $\xi_1 \geq 0$, $\xi_2 \leq 0$, $\xi_3 \leq 0$. In this case, we have $|\xi_1|\geq|\xi_2|$, $|\xi_1|\geq|\xi_3|$,
		$$
		P_6\left(\xi_1\right):=L_1+L_2=\xi_3\left(2 \xi_3^2+3 \xi_1 \xi_3+3 \xi_1^2-1\right)-L_3,
		$$
		and
		$
		|P_6^{\prime}\left(\xi_1\right)|=|3 \xi_3\left(2 \xi_1+\xi_3\right)| \gtrsim |\xi_3\xi_1|
		$.
		Furthermore, one also has
		$
		\left\langle L_3\right\rangle \gtrsim\left\langle H\right\rangle \gtrsim\left\langle\xi_1 \xi_2 \xi_3\right\rangle.
		$
		We notice that, for $\frac{1}{2}<s<\frac{5}{8}$, one has $  -3+2s < \frac{3}{2} -2s$ and it yields that,
		\begin{align*}
			\left\langle\xi_1\right\rangle^{-1-2s}\left\langle\xi_2\right\rangle^{4-2s}\left|\xi_3\right|^{-3+2s} \lesssim
			\left\langle\xi_1\right\rangle^{\frac{3}{2}-2s}\left\langle\xi_2\right\rangle^{\frac{3}{2}-2s}\left|\xi_3\right|^{-3+2s} \lesssim
			\left|\xi_1 \xi_2 \xi_3\right|^{\frac{3}{2}-2s} \lesssim\left\langle\xi_1 \xi_2 \xi_3\right\rangle^{\frac{3}{2}-2s} \lesssim\left\langle L_1\right\rangle^{\frac{3}{2}-2s} .
		\end{align*}
		Therefore, (\ref{3.20}) is bounded by
		$$
		 \frac{1}{\left\langle L_3\right\rangle^{2-2b}} \int \frac{\left\langle\xi_2\right\rangle^{4-2s}}{\left\langle\xi_1\right\rangle^{2s}\left|\xi_3\right|^{3-2s}|\xi_1|} \frac{|P_6^{\prime}\left(\xi_1\right)|}{\left\langle P_6\left(\xi_1\right)\right\rangle^{2 b}} d \xi_1 
			\lesssim
		 \frac{\left\langle L_3\right\rangle^{\frac{3}{2}-2s}}{\left\langle L_3\right\rangle^{2-2b}},
		$$
		which is finite if $\frac{3}{2}-2s \leq 2-2b$.
		
		\item \textbf{Case B6.3.}  If $\left(\xi_1, \xi_2, \xi_3\right)$ lies in \textbf{(c)}, that is, $\xi_1 \leq 0$, $\xi_2 \geq 0$, $\xi_3 \leq 0$. In this case, we have $|\xi_2|\geq|\xi_1|$, $|\xi_2|\geq|\xi_3|$,
		$$
		P_7\left(\xi_1\right):=L_1+L_2=-3 \xi_1 \xi_3\left(\xi_1+\xi_3\right)-L_3,
		$$
		and
		$
		|P_7^{\prime}\left(\xi_1\right)|= 3|\xi_3|(|\xi_1|+|\xi_2|) \gtrsim |\xi_3\xi_2|.
		$
		Additionally, we point out that
		$
		\left\langle L_3\right\rangle \gtrsim\left\langle H\right\rangle  \gtrsim\left\langle\xi_1 \xi_2 \xi_3\right\rangle.
		$
		Thus, (\ref{3.20}) is bounded by
		$$
		 \frac{1}{\left\langle\xi_3\right\rangle^{2-2s}\left\langle L_3\right\rangle^{2-2b}} \int \frac{\left\langle\xi_2\right\rangle^{4-2s}}{\left\langle\xi_1\right\rangle^{2s}\left|\xi_3 \xi_2\right|} \frac{\left|P_7^{\prime}\left(\xi_1\right)\right|}{\left\langle P_7\left(\xi_1\right)\right\rangle^{2 b}} d \xi_1 .
		$$
		For $\frac{1}{2}<s<\frac{5}{8}$, it yields $-3+2s<\frac{3}{2}-2s$ and $\frac{3}{2}-2s>0$. Thus, based on $|\xi_1|,|\xi_2|,|\xi_3|>1$, one has
		\begin{itemize}
			\item \textbf{Case B6.3.1} $|\xi_2| \sim |\xi_1| \gtrsim |\xi_3|$
			\begin{align*}
				\left\langle\xi_1\right\rangle^{-2s}\left\langle\xi_2\right\rangle^{3-2s}\left|\xi_3\right|^{-3+2s} \lesssim
				\left\langle\xi_1\right\rangle^{\frac{3}{2}-2s}\left\langle\xi_2\right\rangle^{\frac{3}{2}-2s}\left|\xi_3\right|^{-3+2s} \lesssim
				\left|\xi_1 \xi_2 \xi_3\right|^{\frac{3}{2}-2s} \lesssim\left\langle\xi_1 \xi_2 \xi_3\right\rangle^{\frac{3}{2}-2s} \lesssim\left\langle L_1\right\rangle^{\frac{3}{2}-2s} ,
			\end{align*}
			\item \textbf{Case B6.3.2} $|\xi_2| \sim |\xi_3| \gg |\xi_1|$
			\begin{align*}
				\left\langle\xi_1\right\rangle^{-2s}\left\langle\xi_2\right\rangle^{3-2s}\left|\xi_3\right|^{-3+2s} \lesssim
				\left\langle\xi_1\right\rangle^{-2s}\lesssim
				\left|\xi_1 \xi_2 \xi_3\right|^{\frac{3}{2}-2s} \lesssim\left\langle\xi_1 \xi_2 \xi_3\right\rangle^{\frac{3}{2}-2s} \lesssim\left\langle L_1\right\rangle^{\frac{3}{2}-2s} ,
			\end{align*}
		\end{itemize}
		which follows that (\ref{3.20}) is finite owing to $\frac{3}{2}-2s \leq 2-2b$.
	\end{itemize}
	
	The proof is now complete.
	\subsection{Proof of Estimate \eqref{tri1}}
	Based on Tao's setting in \cite{tao2001multilinear}, we introduce an $\left[k ; \mathbb{R}\right]$-multiplier to be any function $m: A_k(\mathbb{R}) \rightarrow \mathbb{C}$ and its norm $\|m\|_{\left[k ; \mathbb{R}\right]}$ to be the best constant such that
	$$
	\left|\int_{A_k(\mathbb{R})} m(\xi_1,\tau_1,..., \xi_k, \tau_k) \prod_{j=1}^n h_j(\xi_j,\tau_j)\right| \leq\|m\|_{\left[k ; \mathbb{R}\right]} \prod_{j=1}^n\left\|h_j\right\|_{L^2\left(\mathbb{R}^{2}\right)},
	$$
	holds true for all the test functions in $\mathbb{R}^{2}$ (see \cite{carvajal2019sharp}).
	Similar to the setting in the proof of {\bf Proposition \ref{bi0}}, based on duality and Plancherel theorem (see (3.48) in \cite{carvajal2019sharp}), the estimate \eqref{tri1} is equivalent to
	\begin{align} \label{star}
		\left|\int_{A_4} \frac{|\xi_4| h_1(\xi_1,\tau_1)h_2(\xi_2,\tau_2)h_3\left(\xi_3, \tau_3\right)h_4\left(\xi_4, \tau_4\right)}{\langle \xi_1 \rangle ^ {s}\langle \xi_2 \rangle ^ {s}\langle \xi_3 \rangle ^ {s}\langle \xi_4 \rangle ^{2-s}\langle L_1\rangle^b\langle  L_2\rangle^b\langle L_3\rangle^b\left\langle L_4\right\rangle^{1-b}}  \right|\lesssim \prod_{j=1}^4\left\|h_j\right\|_{L_{\xi, \tau}^2}, \quad \forall h_1,h_2,h_3,h_4 \in \mathscr{S}(\mathbb{R}^2),
	\end{align}
	 where $L_i$, $i=1,2,3,4$, are defined similarly as in {\bf Proposition \ref{bi0}}. Thus, it suffices to show		
	$\|m\|_{\left[4 ; \mathbb{R}\right]} \lesssim 1$ for $b=\frac{1}{2}+\epsilon$ with $\epsilon \in (0, \frac{1}{4})$,
	where,
	\begin{align} \label{m}
		m  \left(\xi_1, \tau_1, \cdots, \xi_4, \tau_4\right) 
		:=\frac{|\xi_4|\left\langle\xi_4\right\rangle^{s-2}}{\left\langle\xi_1\right\rangle^s\left\langle\xi_2\right\rangle^s\left\langle\xi_3\right\rangle^s\left\langle L_1\right\rangle^{\frac{1}{2}+\epsilon}\left\langle L_2\right\rangle^{\frac{1}{2}+\epsilon}\left\langle L_3\right\rangle^{\frac{1}{2}+\epsilon}\left\langle L_4\right\rangle^{\frac{1}{2}- \epsilon}} .
	\end{align}
	Note that $-\xi_4=\xi_1+\xi_2+\xi_3 \Longrightarrow\left\langle\xi_4\right\rangle \leq\left\langle\xi_1\right\rangle+\left\langle\xi_2\right\rangle+\left\langle\xi_3\right\rangle$. Thus, for $s-\frac{1}{2} > 0$, one has
	\begin{align} \label{xi4}
		|\xi_4|\left\langle\xi_4\right\rangle^{s-2} \leq\left\langle\xi_4\right\rangle^{s-1} \lesssim\left\langle\xi_4\right\rangle^{-\frac{1}{2}}\left(\left\langle\xi_1\right\rangle^{s-\frac{1}{2}}+\left\langle\xi_2\right\rangle^{s-\frac{1}{2}}+\left\langle\xi_3\right\rangle^{s-\frac{1}{2}}\right) .
	\end{align}
	From (\ref{m}) and (\ref{xi4}), we get
	$$
	\begin{aligned}
		m \leq  \frac{\left\langle\xi_4\right\rangle^{-\frac{1}{2}}\left\langle\xi_1\right\rangle^{s-\frac{1}{2}}+\left\langle\xi_4\right\rangle^{-\frac{1}{2}}\left\langle\xi_2\right\rangle^{s-\frac{1}{2}}+\left\langle\xi_4\right\rangle^{-\frac{1}{2}}\left\langle\xi_3\right\rangle^{s-\frac{1}{2}}}{\left\langle\xi_1\right\rangle^s\left\langle\xi_2\right\rangle^s\left\langle\xi_3\right\rangle^s\left\langle L_1\right\rangle^{\frac{1}{2}+\epsilon}\left\langle L_2\right\rangle^{\frac{1}{2}+\epsilon}\left\langle L_3\right\rangle^{\frac{1}{2}+\epsilon}\left\langle L_4\right\rangle^{\frac{1}{2}- \epsilon}} 
		=  : J_1+J_2+J_3 
	\end{aligned}
	$$
	where
	\begin{align*}
		 J_1 &=\frac{\left\langle\xi_4\right\rangle^{-\frac{1}{2}}\left\langle\xi_1\right\rangle^{s-\frac{1}{2}}}{\left\langle\xi_1\right\rangle^{\frac{1}{2}}\left\langle\xi_2\right\rangle^s\left\langle\xi_3\right\rangle^s\left\langle\xi_4\right\rangle^{\frac{1}{2}}\left\langle L_1\right\rangle^{\frac{1}{2}+\epsilon}\left\langle L_2\right\rangle^{\frac{1}{2}+\epsilon}\left\langle L_3\right\rangle^{\frac{1}{2}+\epsilon}\left\langle L_4\right\rangle^{\frac{1}{2}- \epsilon}} \\
		& =\left[\left\langle\xi_1\right\rangle^{-\frac{1}{2}}\left\langle L_1\right\rangle^{-\frac{1}{2}-\epsilon}\left\langle\xi_2\right\rangle^{-s}\left\langle L_2\right\rangle^{\frac{1}{2}+\epsilon}\right]\left[ \left\langle\xi_4\right\rangle^{-\frac{1}{2}}\left\langle L_4\right\rangle^{-\frac{1}{2}+ \epsilon}\left\langle\xi_3\right\rangle^{-s}\left\langle L_3\right\rangle^{-\frac{1}{2}-\epsilon} \right]\\
		& \leq n\left(\xi_1, \tau_1, \xi_2, \tau_2\right) n\left(\xi_4, \tau_4, \xi_3, \tau_3\right), 
	\end{align*}
	\begin{align*}
		 J_2 &=\left[\left\langle\xi_4\right\rangle^{-\frac{1}{2}}\left\langle L_4\right\rangle^{-\frac{1}{2}+\epsilon}\left\langle\xi_3\right\rangle^{-s}\left\langle L_3\right\rangle^{-\frac{1}{2}-\epsilon}\right]\left[ \left\langle\xi_2\right\rangle^{-\frac{1}{2}}\left\langle L_2\right\rangle^{-\frac{1}{2}- \epsilon}\left\langle\xi_1\right\rangle^{-s}\left\langle L_1\right\rangle^{-\frac{1}{2}-\epsilon} \right]\\
		& \leq n\left(\xi_4, \tau_4, \xi_3, \tau_3\right) n\left(\xi_2, \tau_2, \xi_1, \tau_1\right), 
	\end{align*}
	\begin{align*}
		J_3
		&=\left[\left\langle\xi_4\right\rangle^{-\frac{1}{2}}\left\langle L_4\right\rangle^{-\frac{1}{2}+\epsilon}\left\langle\xi_2\right\rangle^{-s}\left\langle L_2\right\rangle^{-\frac{1}{2}-\epsilon}\right]\left[ \left\langle\xi_3\right\rangle^{-\frac{1}{2}}\left\langle L_3\right\rangle^{-\frac{1}{2}- \epsilon}\left\langle\xi_1\right\rangle^{-s}\left\langle L_1\right\rangle^{-\frac{1}{2}-\epsilon} \right]\\
		& \leq n\left(\xi_4, \tau_4, \xi_2, \tau_2\right) n\left(\xi_3, \tau_3, \xi_1, \tau_1\right) . 
	\end{align*}
	with
	$$
	\begin{aligned}
		n\left(\xi_i, \tau_i, \xi_j, \tau_j\right)=\frac{1}{\left\langle\xi_i\right\rangle^{\frac{1}{2}}\left\langle L_i\right\rangle^{\frac{1}{2}- \epsilon}\left\langle\xi_j\right\rangle^s\left\langle L_j\right\rangle^{\frac{1}{2}+\epsilon}}, \quad \mbox{for } i,j=1,2,3,4 \mbox{ and }  i\neq j.
	\end{aligned}
	$$
	Hence, using comparison principle, permutation and composition properties (see, respectively, Lemmas 3.1,3.3 and 3.7 in \cite{tao2001multilinear}), it suffices to show $\left\|n\right\|_{\left[3 ; \mathbb{R} \right]} \lesssim 1$, which is equivalent to establish the bilinear estimate
	\begin{align} \label{bi}
		\|u v\|_{L^2\left(\mathbb{R}^2\right)} \lesssim\|u\|_{\frac{1}{2}, \frac{1}{2}- \epsilon}\|v\|_{X_{s, \frac{1}{2}+\epsilon}},
	\end{align}
	with $s>\frac{1}{2}$ and $0<\epsilon<\frac{1}{4}$.
	\subsubsection{Proof of (\ref{bi})}
	Similar to the setting in the proof of \textbf{Proposition \ref{bi0}}, based on duality and Plancherel theorem, (\ref{bi}) is equivalent to
	\begin{align} \label{bi1}
		\int_{A_3} \left|\frac{h_1(\xi_1,\tau_1)h_2(\xi_2,\tau_2)h_3(\xi_3,\tau_3)}{\langle \xi_1 \rangle ^ {\frac{1}{2}}\langle \xi_2 \rangle^{s}\langle L_1\rangle^{\frac{1}{2}- \epsilon}\langle  L_2\rangle^{\frac{1}{2}+\epsilon}}  \right|\lesssim \prod_{j=1}^3\left\|h_j\right\|_{L_{\xi, \tau}^2}. 
	\end{align}
	Then, it  suffices to consider Cases \textbf{(a)}, \textbf{(b)} and \textbf{(c)}  lying in sets   \textbf{A} and \textbf{B}. 
	
	Firstly, for Case \textbf{A}, we can write
	$$
	L_1=\tau_1-\left|\xi_1\right|^3+\frac{1}{2}\left|\xi_1\right| , \quad L_2=\tau_2-\left|\xi_2\right|^3+\frac{1}{2}\left|\xi_2\right| .
	$$
	\textbf{Case A1.} $\left|\xi_3\right| \leq 1$. In this case, one has
	\begin{align*}
		\text { LHS of }(\ref{bi1})  \lesssim \iint \left|h_1\right|\left(\iint \frac{d \xi_3 d \tau_3}{\langle \xi_1 \rangle \langle \xi_2 \rangle^{2s}\langle L_1\rangle^{1- 2\epsilon}\langle  L_2\rangle^{1+2\epsilon}}\right)^{\frac{1}{2}}\left(\iint h_2^2 h_3^2 d \xi_3 d \tau_3\right)^{\frac{1}{2}} d \xi_1 d \tau_1 .
	\end{align*}
	According to Lemma \ref{lem1}, for $s>\frac{1}{2}$, $1-2\epsilon>0$ and $|\xi_3|\leq 1$, one has
	\begin{align*}
			\iint \frac{d \xi_3 d \tau_3}{\langle \xi_1 \rangle \langle \xi_2 \rangle^{2s}\langle L_1\rangle^{1- 2\epsilon}\langle  L_2\rangle^{1+2\epsilon}} 
		\lesssim 
		\int \frac{d\xi_3  d\tau_3}{\langle L_1\rangle^{1- 2\epsilon}\langle  L_2\rangle^{1+2\epsilon}} 
		\lesssim
		\int_{|\xi_3|\leq 1} \frac{d \xi_3}{\left\langle L_1+L_2\right\rangle^{1- 2\epsilon}} \lesssim 1 .
	\end{align*}
	Hence, (\ref{bi1}) is achieved through Cauchy-Schwarz inequality,
	\begin{align*}
		\text { LHS of }(\ref{bi1}) & \lesssim \iint\left|h_1\right|\left(\iint h_2^2 h_3^2 d \xi_3 d \tau_3\right)^{\frac{1}{2}} d \xi_1 d \tau_1 
		\lesssim \prod_{j=1}^3\left\|h_j\right\|_{L_{\xi,\tau}^2}.
	\end{align*}
	\textbf{Case A2.} $\left|\xi_3\right| > 1$. In this case, we notice that
	\begin{align*}
		\text { LHS of }(\ref{bi1}) \lesssim \iint \left|h_3\right|\left(\iint \frac{d \xi_1 d \tau_1}{\langle \xi_1 \rangle \langle \xi_2 \rangle^{2s}\langle L_1\rangle^{1- 2\epsilon}\langle  L_2\rangle^{1+2\epsilon}}\right)^{\frac{1}{2}}\left(\iint h_1^2 h_2^2 d \xi_1 d \tau_1\right)^{\frac{1}{2}} d \xi_3 d \tau_3 .
	\end{align*}
	According to Lemma \ref{lem1} , for $0<\epsilon<\frac{1}{4}$ and $s>\frac{1}{2}$, one has
	
\begin{equation}\label{tri}
		\iint \frac{d \xi_1 d \tau_1}{\langle \xi_1 \rangle \langle \xi_2 \rangle^{2s}\langle L_1\rangle^{1- 2\epsilon}\langle  L_2\rangle^{1+2\epsilon}} 
	\lesssim
	 \int \frac{d \xi_1}{\langle \xi_1 \rangle \langle \xi_2 \rangle^{2s}\left\langle L_1+L_2\right\rangle^{1- 2\epsilon}}   
	 \lesssim
	  \int \frac{d \xi_1}{\left\langle L_1+L_2\right\rangle^{1- 2\epsilon}}.  	
\end{equation}
	\begin{itemize}
		\item \textbf{Case A2.1.} If $\left(\xi_1, \xi_2, \xi_3\right)$ lies in \textbf{(a)}, that is, $\xi_1 \geq 0$, $\xi_2 \geq 0$, $\xi_3 \leq 0$. In this case, we have 
		$$
		\begin{aligned}
			Q_1\left(\xi_1\right)  :=L_1+L_2
			=3 \xi_3 \xi_1^2+3 \xi_3^2 \xi_1-L_3.
		\end{aligned}
		$$
		According to Lemma \ref{lem3}, for $1- 2\epsilon > \frac{1}{2}$ and $|\xi_3|>1$, one has,
		$$
		\int \frac{d \xi_1}{\left\langle L_1+L_2\right\rangle^{1- 2\epsilon}}\lesssim \int \frac{1}{\left\langle Q_1\left(\xi_1\right)\right\rangle^{1- 2\epsilon}} d \xi_1
		\lesssim
		\frac{1}{\left| \xi_3\right|^{\frac{1}{2}}}
		\lesssim 1.
		$$
		\item \textbf{Case A2.2.}  If $\left(\xi_1, \xi_2, \xi_3\right)$ lies in \textbf{(b)} and \textbf{(c)},  one can proceed the proof as in \textbf{Case A2.1.} to show \eqref{tri} is bounded.
	\end{itemize}
	Next, we proceed to examine the case that $\left(\tau_1, \tau_2, \tau_3\right)$ lies in region $\mathbf{B}$.
\\
	\textbf{Case B1.} If $\left|\xi_3\right| \leq 1$, it is similar to \textbf{Case A1}, therefore omitted. 
	\\
	\textbf{Case B2.} If $\left|\xi_3\right|>1$, it suffices to show the bound for
	\begin{align}\label{b2}
		\sup _{\xi_3, \tau_3} \int \frac{d \xi_1}{\left\langle L_1+L_2\right\rangle^{1- 2\epsilon}},
	\end{align}
	
	\begin{itemize}
		\item \textbf{Case B2.1.}  If $\left(\xi_1, \xi_2, \xi_3\right)$ lies in \textbf{(a)}, it is similar to \textbf{Case A2.2}.
		\item \textbf{Case B2.2.}  If $\left(\xi_1, \xi_2, \xi_3\right)$ lies in \textbf{(b)} or \textbf{(c)}, this case can be addressed same as \textbf{Case A2.1}. 
	\end{itemize}
	The proof is now complete. 
	
	\smallskip
	
	\centerline{\textbf{Acknowledgements}}

	S. Li is supported by the National Natural Science Foundation of China (no. 12001084 and no. 12071061).


\end{document}